\documentclass[reqno, 10pt]{amsart}  
\usepackage{amsmath}
\usepackage{amssymb}
\usepackage{hyperref}
\hypersetup{colorlinks=true}

\newtheorem{prop}{Proposition}[section]
\newtheorem{lemma}[prop]{Lemma}

\newtheorem{cor}[prop]{Corollary}

\newtheorem{definition}[prop]{Definition}

\newtheorem*{main_thm_sdp}{Main Theorem for Semidefinite Programming}
\newtheorem*{main_thm_hyperbolic}{Main Theorem for Hyperbolic Programming}
\newtheorem*{garding_thm}{G\aa rding's Theorem}
\newtheorem*{helton_vinnikov_thm}{Helton-Vinnikov Theorem}

\newcommand{\shrink}[1]{ {\scriptstyle {\textstyle {#1} } } }
\newcommand{\smfrac}[2]{ \shrink{ \frac{#1}{#2} } }

\newcommand{\lin}{\langle}
\newcommand{\rin}{\rangle}
\newcommand{\Sym}{\mathbb{S}^n}
\newcommand{\Symp}{\Sym_{\plus}} 
\newcommand{\Sympp}{\Sym_{\pplus}}

\newcommand{\diag}{\textrm{diag}}

\newcommand{\plus}{{\scriptscriptstyle +}}
\newcommand{\pplus}{\plus \plus} 
\newcommand{\Lambdap}{\Lambda_{\plus}}
\newcommand{\Lambdapp}{\Lambda_{\plus \plus}}

\newcommand{\tr}{\mathrm{ {\bf  tr}}} 
\newcommand{\sdp}{\mathrm{SDP}}

\newcommand{\hp}{\mathrm{HP}}

\newcommand{\swath}{\mathrm{Swath}} 
 
\renewcommand{\int}{\mathrm{int}}

\newcommand{\qp}{\mathrm{QP}}

\begin{document}

\title[A Primal Affine-Scaling Algorithm]{A Polynomial-Time Affine-Scaling Method \\ for Semidefinite and Hyperbolic Programming}
\author{James Renegar and Mutiara Sondjaja}
\address{School of Operations Research and Information Engineering, Cornell University}
\address{Courant Institute of Mathematical Sciences, New York University}

\begin{abstract}
We develop a natural variant of Dikin's affine-scaling method, first for semidefinite programming and then for hyperbolic programming in general.  We match the best complexity bounds known for interior-point methods.

All previous polynomial-time affine-scaling algorithms have been for conic optimization problems in which the underlying cone is symmetric.  Hyperbolicity cones, however, need not be symmetric.  Our algorithm is the first polynomial-time affine-scaling method not relying on symmetry.
\end{abstract}

\maketitle

\vspace{-8mm}

\section{{\bf Introduction}} \label{sect.a}

Dikin's affine-scaling method for linear programming (\cite{dikin1967iterative}) received considerable attention during the 1980's as a natural variant of Karmarkar's projective-scaling algorithm (\cite{Karmarkar:1984jb}).  Compared with Karmarkar's algorithm, the affine-scaling method is conceptually simple.  Indeed, for a linear program
\[ \left. \begin{array}{rl}
   \min & c^T x \\
   \textrm{s.t.} & Ax = b \\
  & x \geq 0  \end{array} \quad \right\} \, \textrm{LP} \; , \]
the algorithm can be described as follows:
\begin{quote}
  Given a current iterate $ e $ which is assumed to be strictly feasible, perform the ``scaling'' $ x \mapsto  z := E^{-1} x $ where $ E = \textrm{diag($e$)} $,  thereby obtaining the equivalent linear program 
  \[ \begin{array}{rl} \min_z & (E \, c)^T z \\
                 \textrm{s.t.} & (A \, E) \, z = b \\
                 & z \geq 0 \; , \end{array} \]
for which the current iterate is $ \mathbf{1} $, the vector of all ones.  Letting $ d $ denote the Euclidean projection of the objective vector $ E c $  onto the nullspace of the constraint matrix $ A E $, step from $ {\bf 1} $ to $ {\bf 1} - \smfrac{1}{\|d\|}   d $  (a step of unit length).  Finally, reverse the scaling to obtain $ \bar{e}  := E \, \left(   {\bf 1} - \smfrac{1}{\|d\|} d \right)  $, and define this to be the iterate subsequent to the current iterate, $ e $. (Clearly, $ \bar{e} $ is feasible. Moreover, it happens that if some coordinate of $ \bar{e} $ is zero, then $ \bar{e} $ is optimal, and the algorithm terminates.)   
\end{quote} \vspace{1mm}

An enlightening (and equivalent) alternative description of Dikin's method relies on the inner products $ \lin u, v \rin_e := \sum \frac{u_j  v_j}{e_j^2} $ defined for all $ e $ with strictly-positive coordinates.  Letting $ \| v \|_e := \lin v, v \rin_e^{1/2} \; , $ the induced norm, here is the alternative description:
\begin{quote}
  Given a current iterate $ e $ which is assumed to be strictly feasible, let $ \bar{e} $ (the next iterate) be the optimal solution to 
  \[ \begin{array}{rl}
 \min_x & c^T x \\
      \textrm{s.t.} & Ax = b \\
        & x \in \bar{B}_e(e,1)  \; ,  
        \end{array} \]
where $ \bar{B}_e( e, 1) = \{ x: \|x - e \|_e \leq 1 \} $, the $ \| \;\; \|_e $-unit ball centered at $ e $ -- often referred to as the ``Dikin ellipsoid'' centered at $ e $.
\end{quote}

Whereas the initial description of the algorithm makes evident the origin of the terminology ``affine-\underline{scaling} method,'' the alternative description highlights the simple idea underlying the algorithm: Replace the non-negativity constraints $ x \geq 0 $ in LP with the (computationally easier) restriction that $ x $ lie in the Dikin ellipsoid at $ e $, the largest ellipsoid centered at $ e $ that both is contained in the non-negative orthant and whose axes are parallel to the coordinate axes. \vspace{1mm}

Dikin introduced the algorithm in the mid-1960's but it remained unknown in the West until after the framework had been rediscovered in the mid-1980's by Vanderbei, Meketon and Freedman \cite{Vanderbei:1986gy}, who rather than implementing a step of unit length, chose to step from $ \mathbf{1} $ to $ \mathbf{1} - \frac{\gamma}{ \max_i d_i} \,  d $ for fixed $ 0 < \gamma < 1 $, that is, chose to step a fraction of the distance, in direction $ -d $, to the boundary of the feasible region.  They showed this step choice results in a sequence of iterates that converges to an optimal solution under the assumption that LP is both primal and dual nondegenerate, whereas Dikin  proved his choice leads to optimality assuming only dual nondegeneracy \cite{dikin1974speed, Vanderbei:1990uj}.  \vspace{1mm}

For linear programming, the most definitive results regarding step choice are due to Tsuchiya and Muramatsu \cite{TSUCHIYA:1995cb}, who showed that choosing $ 0 < \gamma \leq  2/3 $ always results in convergence to optimality regardless of whether nondegeneracy holds.  Moreover, in a sense that we do not describe here, they showed the bound $ \gamma \leq 2/3 $ in general is the best possible. (Also see \cite{monteiro1993simplified}.)  \vspace{1mm}

Despite the appeal of its simplicity, and despite the various convergence results noted above, the affine-scaling method has long been believed {\em not} to provide a polynomial-time algorithm for linear programming.  The first rigorous findings supporting this view were developed in the mid-1980's by Megiddo and Shub \cite{Megiddo:1989wq}, who displayed for a Klee-Minty cube that given $ \epsilon > 0 $, by beginning at a particular interior point dependent on $ \epsilon $, if the affine-scaling steps are made infinitesimally -- thereby resulting in a smooth path -- the path comes within distance $ \epsilon $ of each of the (exponentially many) vertices of the cube. \vspace{1mm}

There have been, however, a number of interior-point methods which share spirit with Dikin's affine-scaling method and for which polynomial-time bounds have been established, for linear programming (Monteiro, Adler and Resende \cite{Monteiro:1990vn}; Jansen, Roos and Terlaky \cite{Jansen:1996bb}; Sturm and Zhang \cite{Sturm:1996va}), then for semidefinite programming (de Klerk, Roos and Terlaky \cite{deKlerk:1998dc}; Berkelaar, Sturm and Zhang \cite{Berkelaar:1999tz}), and later for symmetric-cone programming in general (Chua \cite{Chua:2007cu}). The complexity bounds in  \cite{Sturm:1996va}, \cite{Berkelaar:1999tz} and \cite{Chua:2007cu} match the best available for interior-point methods. \vspace{1mm}

All of these algorithms share spirit with Dikin's method in that a step is determined by solving the optimization problem in which the conic constraint has been replaced by an ellipsoidal constraint or, for \cite{Sturm:1996va}, \cite{Berkelaar:1999tz} and \cite{Chua:2007cu}, by an ellipsoidal cone constraint (i.e., a cone $ \{ tx: t \geq 0 \textrm{ and } x \in {\mathcal E} \} $ where $ {\mathcal E} $ is an ellipsoid).  Each of these algorithms is reasonably called an affine scaling method.    \vspace{1mm}  

(The idea of using ellipsoidal cones rather than ellipsoids in designing interior-point methods seems first to have appeared in technical reports of Megiddo \cite{megiddo1985variation} and Todd \cite{todd1991another}, albeit in relation to potential reduction methods.)  
\vspace{1mm}

All of the above polynomial-time affine-scaling methods are primal-dual algorithms, with complexity analysis depending heavily on properties of symmetric cones (a.k.a. self-scaled cones \cite{nesterov1997self}).  With the exception of the algorithm in \cite{Chua:2007cu}, even motivation of the computational steps rests on properties only known to be possessed by symmetric cones.  (The algorithm in \cite{Chua:2007cu} is unusual also in that it alternates between the primal and dual problems, rather than merging the problems into a single problem in ``primal-dual space.'') \vspace{1mm} 

We do not provide details of these polynomial-time affine-scaling methods because our primary goal is to devise an efficient and easily motivated affine-scaling method which applies to optimization problems defined over a far broader class of cones than just symmetric ones.  Nonetheless, we do remark that our work very much depends on perspectives first put forward by Chua \cite{Chua:2007cu} in developing his algorithm. \vspace{1mm} 

We present and analyze a primal algorithm -- or it could be cast as a dual algorithm -- which matches the best complexity bounds of interior-point methods.  When, for example, the underlying cone is $ \mathbb{R}^n_+ $ (non-negative orthant) or $ \Symp $ (positive semidefinite cone), the complexity bound is $ O( \sqrt{n}) $ iterations to ``halve the duality gap'' (the algorithm puts emphasis on primal iterates but also provides easily-computable feasible dual iterates, and thus naturally gives duality gaps with which to measure progress). \vspace{1mm}

Our algorithm shares spirit with Dikin's method in that at each iteration, the cone in the primal optimization problem is replaced by a simpler set -- in our case the set is an ellipsoidal cone centered at the current (primal) iterate, whereas for Dikin it was an ellipsoid.  Our algorithm differs in spirit, however, in that at each iteration, the ellipsoidal cone is chosen to contain the original primal cone rather than be contained in it.  In other words, at each iteration the original primal problem is ``relaxed'' to an ellipsoidal-cone optimization problem with larger feasible region.  A consequence is that at each iteration, the original dual optimization problem is ``restricted'' to an ellipsoidal-cone optimization problem with smaller feasible region. \vspace{1mm}

The ellipsoidal-cone optimization problems are easily solved simultaneously (when they have optimal solutions), as the first-order optimality conditions for the primal yield the dual optimal solution.  Moreover, the optimal solution of the dual ellipsoidal-cone problem provides a feasible iterate for the original dual problem, simply because the dual ellipsoidal-cone problem is a restriction of the original dual problem.  By contrast, the optimal solution of the primal ellipsoidal-cone problem is infeasible for the original primal problem, and hence is not used as the next primal iterate. \vspace{1mm}

Instead, the next primal iterate is obtained as a convex combination of the current primal iterate -- at which the primal ellipsoidal cone is centered -- and the optimal solution of the primal ellipsoidal-cone problem (this optimal solution lies on the boundary of the primal ellipsoidal cone).  The specific weights in the convex combination are determined in an elegant computation depending also on the optimal solution to the dual ellipsoidal-cone problem.  In this way, the dual iterates help guide the choice of the primal iterates, and there can result long steps (i.e., this is not a short-step method, although short steps sometimes occur). \vspace{1mm}

We first develop the algorithm and analysis for semidefinite programming (sections \ref{sect.b}--\ref{sect.f}).  Although the algorithm is easily motivated far beyond semidefinite programming, not so for the choice of step length, and it is unclear the extent to which the complexity analysis can be extended.  However, we are able to leverage the semidefinite programming step-length selection, and complexity analysis, so as to apply to all hyperbolic programs. (An overview of hyperbolic programming is provided in section \ref{sect.g} for readers unacquainted with the topic.)  A key piece of the leveraging is done using a deep theorem whose statement is easily understood, the Helton-Vinnikov Theorem (discussed in section \ref{sect.j}).   \vspace{1mm}

Hyperbolicity cones -- the cones underlying hyperbolic programming -- are in general not self-dual (let alone symmetric), as was observed by G\"{u}ler \cite{guler1997hyperbolic}  when he introduced hyperbolic programming.  Thus, the fact that the algorithm and complexity results pertain to all hyperbolic programs makes evident that by no means is this a primal-dual algorithm, even if the dual iterates help guide the choice of primal iterates.  However, our interest in having the results apply to all of hyperbolic programming goes beyond the academic desire of making rigorous that the algorithm is primal, or dual, but not primal-dual (i.e., not, from some perspective, symmetric). \vspace{1mm}

A key feature of hyperbolic programming is that for any hyperbolic program -- say, a semidefinite program -- there is a natural hierarchy of hyperbolic program relaxations of it, with each step in the hierarchy resulting in an optimization problem whose boundary structure is somehow simpler than that of the previous ones (see \cite{renegar2006hyperbolic}).  In particular, with each step in the hierarchy for, say, a semidefinite program, the parameter controlling the efficiency of interior-point methods is lessened by one, and thus worst-case complexity is reduced from $ O(\sqrt{n}) $ iterations to halve the duality gap, to $ O(\sqrt{n-1}) $ iterations to halve the gap, all the way down to $ O(\sqrt{2}) $ iterations.  If the optimal solution for a relaxation happens to be nearly positive semidefinite, then it serves as an approximate optimal solution to the original semidefinite program, albeit one obtained with relatively small cost. \vspace{1mm}

The point here is that hyperbolic programming, besides possessing interesting theoretical characteristics (such as not supporting primal-dual algorithms), has characteristics that might prove relevant to approximately solving, say, much larger general semidefinite programs than can currently be handled.   \vspace{1mm}

In the following section -- \S\ref{sect.b}  -- we begin motivating the semidefinite programming algorithm from a primal perspective, and present results which underscore the reasonableness of the basic ideas behind the algorithm.  However, this development can go only so far, and it soon becomes apparent that additional ideas are needed, especially for a complexity analysis.  In section \ref{sect.c}  we are thus led to bring duality into the picture, even though the algorithm is primal.  In particular, we make evident how the duality gap at each iteration has meaning even from a primal-only perspective, and we state the \hyperlink{targ_main_thm_sdp}{Main Theorem for Semidefinite Programming}, which includes a specification of the computation made at each iteration to determine the step size.  However, it is not until the section thereafter -- \S\ref{sect.d}  -- that we go deep enough to begin revealing the geometry underlying the choice of step size.  Exploration of this geometry and performing the complexity analysis become intertwined until the proof of the main theorem is completed in section \ref{sect.f}.  \vspace{1mm}

Then, as already indicated, we turn attention to hyperbolic programming in general, providing a motivational overview for the reader, followed by the leveraging of our earlier results so as to extend to all hyperbolic programs. The \hyperlink{targ_main_thm_hp}{Main Theorem for Hyperbolic Programming} is stated in section \ref{sect.i} and proven in section \ref{sect.k}.  \vspace{1mm}

In closing the introduction, we wish draw attention to work of Litvinchev that we only recently happened upon, but which bears important similarities, as well as differences, with our approach.  In particular, in \cite{Litvinchev:2003iy}  he develops a potential-reduction method for linear programming in which at each iteration, the primal problem is relaxed to an ellipsoidal-cone problem.  If the relaxation has an optimal solution, the current iterate is moved towards the optimal solution, just as in our algorithm.  If no optimal solution exists for the relaxation, the current iterate is moved towards a feasible point for the relaxation, a feasible point whose objective value is at least as good as a lower bound on the linear program's objective value (the lower bound is updated every iteration). In either case, the step length from the current iterate is determined by minimizing a potential function.  \vspace{1mm}

(While the reports of Megiddo \cite{megiddo1985variation} and Todd \cite{todd1991another} relied on ellipsoidal cones in investigating potential-reduction, their ellipsoidal cones were contained in the non-negative orthant, whereas Litvinchev's cones, like ours, circumscribe \vspace{1mm} the orthant.)

Litvinchev does not develop any duality perspectives in presenting and analyzing his algorithm.  In particular, a dual iterate is not specified, and thus neither is a duality gap.  Instead of establishing an iterative result such as ``the duality gap halves in $ O( \sqrt{n}) $ iterations,'' he establishes the asymptotic result that $ O( n \, \log(1/ \epsilon )) $ iterations suffice to obtain a feasible point with objective value within $ \epsilon $ of optimality, that is, a result similar to the one proven by Karmarkar \cite{Karmarkar:1984jb}.  (In this kind of result, the big-$ O $ necessarily depends on the particular instance (c.f., \cite[Cor.~3.1]{todd1997potential}). By contrast, big-$ O $'s for halving the duality gap can be universal.  We precisely specify the constant for our big-$ O $, making clear it is universal and reasonable in size.) \vspace{1mm}

Although our results are far more general, it seems entirely plausible that the results of Litvinchev could be greatly extended by making use of general potential functions (c.f., \cite[Chap.~4]{nesterov1994interior}).  However, it seems unlikely that his ``$ n $'' could be reduced to ``$ \sqrt{n} $'', and unlikely that there is a rich underlying duality theory in which dual iterates can be seen as helping determine the choice of primal iterates.  Nonetheless, Litvinchev's results are accomplished in significantly fewer pages than our results.  His work deserves to be highlighted here due its primal focus, and due to its use of the same ellipsoidal-cone relaxations.  \newpage

\section{{\bf  Preliminaries and Motivation \\ in the Case of Semidefinite Programming}}  \label{sect.b}

For $ C, A_1, \ldots, A_m \in \Sym $ ($ n \times n $ symmetric matrices), and $ b \in \mathbb{R}^{m} $, consider the semidefinite program   
\[ \left. \begin{array}{rl}
\min_{X \in \Sym} & \tr(CX) \\
\textrm{s.t.} & {\mathcal A}  (X) = b \\
  & X \succeq 0 \end{array} \right\} \sdp \]
where $ \tr $ is the trace operator, where  $ {\mathcal A} (X) :=  ( \tr(A_1 X), \ldots, \tr(A_m X)) $, and where $ X \succeq 0 $ is shorthand for $ X \in \Symp $ (cone of positive semidefinite matrices).  The dual problem is
\[ \left. \begin{array}{rl}
   \max_{(y,S)} & b^T y \\ \textrm{s.t.} & \sum_i y_i A_i + S = C \\
      & S \succeq 0 \end{array} \right\} \, \sdp^* \]  
      Sometimes we write $ {\mathcal A}^* y $ instead of $ \sum_i y_i A_i $. \vspace{1mm}

As is standard in the interior-point method literature, we assume both problems are strictly feasible, that is, there exist feasible points $ X $ and $ (y,S) $ for which $ X $ and $ S $ are strictly positive definite.  Then, as is well known, both problems have optimal solutions and both problems have the same optimal value.  Let opt\_val denote the optimal value.       \vspace{1mm}

We assume $ b \neq 0 $ (otherwise $ X = 0 $ is optimal for SDP), and we assume $ C $ is not a linear combination of $ A_1, \ldots, A_m $ (otherwise all points feasible for SDP are optimal).  The latter assumption implies, of course, that no strictly feasible point for SDP is optimal.  \vspace{1mm}

To each $ E \in \Sympp $ (strictly positive-definite) there is associated an inner product on $ \Sym $, given  by 
\begin{align*}
   \lin U, V \rin_E & := \tr(U E^{-1}   V E^{-1}) \\
   & = \tr( E^{-1} U E^{-1} V ) \\ 
     & = \tr( \, (E^{-1/2} U E^{-1/2}) \, (E^{-1/2} V E^{-1/2})  \, ) \; ,
     \end{align*} 
     where $ E^{1/2} $ is the strictly positive-definite matrix satisfying $ E = E^{1/2} E^{1/2} $.
As with every inner product, there is an induced norm:  $ \big\| V \big\|_E = \lin V, V \rin_E^{1/2} $.  \vspace{1mm}

For positive value $ \alpha $, consider the ellipsoidal cone -- or convex ``quadratic'' cone -- specified by
\begin{align}
     K_E( \alpha ) & := \{ X:  \lin E, X \rin_E \geq \alpha \| X \|_E \} \label{eqn.ba}
     \\
   & = \{ X: \tr(E^{-1} X) \geq \alpha \, \tr( \, (E^{-1} X)^2 \, )^{1/2} \} \nonumber \\ 
    & =    \{ X: \tr(E^{-1}X)^2 -   \alpha^2 \,  \tr\big(\, (E^{-1/2}XE^{-1/2})^2  \, \big) \textrm{ and } \tr( E^{-1}X ) \geq 0    \} \; .  \nonumber
   \end{align}     
It is seen from (\ref{eqn.ba})   that from the viewpoint of the inner product $ \langle \; , \; \rangle_E $, the cone is circular, consisting precisely of the ``vectors'' $ X  \in \Sym $ for which the angle with $ E $ does not exceed $ \arccos( \alpha / \sqrt{n} )  $  (using that $ \|E\|_E = \sqrt{n} \, $). The set $ K_e( \alpha) $  is nonempty precisely when $ 0 \leq \alpha \leq \sqrt{n} $, and is a regular cone (i.e., pointed with nonempty interior) precisely when $ 0 < \alpha < \sqrt{n} $.  
 \vspace{1mm}
 
 \begin{lemma}  \label{lem.baa}
 If $ E \in \Sympp $ and $ 0 < \alpha \leq 1 $, then $ \Symp \subseteq K_E( \alpha ) \; . $ 
 \end{lemma}
\noindent {\bf Proof:}   
Indeed, assuming $ X \in \Sym$ and letting $ \lambda $ denote the vector of eigenvalues for $ E^{-1} X $ -- the same eigenvalues as for the symmetric matrix $ E^{-1/2} X E^{-1/2} $ -- then
 \[  \lin E, X \rin_E \geq \alpha \, \| X \|_E \quad \Leftrightarrow \quad \sum_j \lambda_j \geq \alpha \, \| \lambda \| \; , \]
which clearly holds if $ \alpha \leq 1 $ and the eigenvalues are non-negative, as they are when $ X $ -- and hence $ E^{-1/2} X E^{-1/2} $ -- is positive semidefinite. \hfill $ \Box $
 \vspace{3mm}
 
 Thus, given $ E \in \Sympp $ and $ 0 < \alpha \leq 1 $, the following quadratic-cone optimization problem is a relaxation of SDP (i.e., it has the same objective vector as SDP and its feasible region contains the feasible region for SDP):  
\[ \left. \begin{array}{rl}
\min_{X} & \tr(CX) \\
\textrm{s.t.} & {\mathcal A} (X) = b \\
  & X \in K_E( \alpha ) \end{array} \right\} \textrm{QP}_E( \alpha)  \]
A benefit is that this problem is easily solved when it has an optimal solution.  Indeed, as any optimal solution must lie in the boundary $ \partial K_E( \alpha ) $ of $ K_E( \alpha ) $ (because, by assumption, $ C $ is not a linear combination of $ A_1, \ldots, A_m $),  necessary first-order conditions for $ X $ to be optimal are
\begin{align}
  & 0 = \overbrace{\smfrac{1}{2} \Big( \tr(E^{-1}X)^2 - \alpha^2 \, \tr\big( \,  (E^{-1}X)^2 \, \big) \Big)}^{\textrm{denote this function by $ X \mapsto   g(X)$}}  \nonumber \\ & \qquad  \qquad    \textrm{(if $ X \in \partial K_E( \alpha) $, then $ X $ satisfies this equation)} \nonumber \\  & {\mathcal A} (X) = b \nonumber \\
   & 0 = \lambda C + {\mathcal A}^* y  +  \overbrace{\tr(E^{-1}X) \, E^{-1}  - \alpha^2 E^{-1} X E^{-1}}^{\textrm{the gradient of $ g $ at $ X $}}  \label{eqn.bb} \\ &
\qquad  \qquad    \textrm{for some $ y \in \mathbb{R}^{m} $ and $ \lambda \in \mathbb{R} $} \; . \nonumber \end{align}
Excluding the first equation gives a system which is linear in the variables $ (X, y, \lambda) $, and which has a solution set of dimension one (assuming,  as we may, that $ A_1, \ldots , A_m $ are linearly independent in the vector space $ \Sym $).  The first equation can then be used to isolate two candidates $ X', X'' $ -- simply substitute a linear parameterization of the one-dimensional set into the first equation and use the quadratic formula to find the roots.  Finally, for those two candidates, checking feasibility  and comparing objective values reveals which one is the optimal solution  -- {\em  \underline{if}}     there is an optimal solution. \vspace{1mm}

However, $ \textrm{QP}_E( \alpha ) $ might not have an optimal solution even if SDP does have one.  Indeed, $ \textrm{QP}_E( \alpha ) $, being a relaxation of SDP (assuming $ 0 \leq \alpha \leq 1 $), might have unbounded optimal value instead.  This leads us to make a definition.
\begin{definition} \label{def.bab} 
 \[  \swath( \alpha) := \{ E \in \Sympp: {\mathcal A} (E) = b \textrm{ and } \textrm{QP}_E( \alpha) \textrm{ has an optimal solution} \} \] 
 \end{definition}
 \vspace{2mm}
 
 These ``swaths'' become broader as $ \alpha $ grows, that is, $ \swath( \alpha_1 ) \subseteq \swath( \alpha_2 ) $ for $ 0 \leq \alpha_1 \leq \alpha_2 \leq \sqrt{n} $.  Indeed,  $ K_E( \alpha_1  ) \supseteq K_E( \alpha_2 ) $, and thus $ \qp_E( \alpha_1 ) $ is a relaxation of $ \qp_E( \alpha_2 ) $.  Hence, if $ \qp_E( \alpha_1 ) $ has an optimal solution, then so does $ \qp_E( \alpha_2 ) $ so long as $ \qp_E( \alpha_2 ) $  is feasible.  But of course $ \qp_E( \alpha_2 ) $ is feasible -- indeed, $ E $ is a feasible point.  \vspace{1mm}

For $ \alpha $ decreasing to zero, the limit of the swaths is most ubiquitous of sets in the interior-point method literature, namely, the central path. \vspace{1mm}

\begin{prop}  \label{prop.bb}
$ \mathrm{Central \,  Path} = \swath(0) $ 
\end{prop}
\noindent {\bf Proof:}  This is Theorem 2 in \cite{renegar2013central}. \hfill $ \Box $
\vspace{3mm}

Proposition~\ref{prop.bb} is purely motivational and plays no role in what follows (otherwise we would include a proof).  However, the proposition provides the reader with the insight that when it is assumed $ E \in \swath( \alpha ) $, the value of $ \alpha $ reflects a bound on the ``proximity'' of $ E $ to the central path.  
\vspace{1mm} 

Here is a setting of interest to us:  Fix $ 0 < \alpha \leq  1 $. Assume we know $ E \in \swath( \alpha ) $, and we want $ E' \in \swath( \alpha ) $ with ``notably better'' objective value than $ E $.  Is there a straightforward manner in which to compute such a point $ E' \,  $? \vspace{1mm}

Well, we have nothing to work with other than that $ \textrm{QP}_E( \alpha ) $ has an optimal solution.  Let's denote the optimal solution by $ X_E $ ($ = X_E( \alpha ) $). As $ 0 < \alpha \leq 1 $, observe that $ \textrm{opt\_val} $ (the optimal value of SDP) satisfies 
\[    \tr(C X_E) \leq  \textrm{opt\_val} <  \tr(C E) \; , \]
the first inequality due to $ \qp_E( \alpha ) $ being a relaxation of SDP, and the second inequality due to $ E $ being strictly feasible for SDP.    Thus, to obtain $ E' $ with better objective value than $ E $, we can choose a convex combination of $ E $ and $ X_E $, that is, choose $ E' = E(t) := \smfrac{1}{1+t} \big( E + t X_E \big) $ for some positive $ t  \; . $ Then, in fact, we will have
\begin{equation}  \label{eqn.bc}
  \frac{\tr(C \, E(t) ) - \textrm{opt\_val} }{\tr(C \, E ) - \textrm{opt\_val}  }  \leq \frac{ \tr( C E(t)) - \tr(C X_E)   }{\tr(CE) - \tr(CX_E)  } = \frac{1}{1+t}  
\end{equation} 
(the inequality due to $ \frac{s_1 - s_3 }{ s_2 - s_3 } $ being a decreasing function of $ s_3 $ when $ s_1, s_3  < s_2 $).  \vspace{1mm}

However,  besides $ E' $ having better objective value than $ E $, we also want $ E' $ to lie in $ \swath( \alpha ) $.   The stage is set for the following result.  \vspace{1mm}

\begin{prop} \label{prop.bc}
  Assume $ 0 < \alpha \leq  1 $ and $ E \in \swath( \alpha ) $.  Let $ X_E = X_E( \alpha ) $ be the optimal solution of $ \qp_E( \alpha ) $.  If $ t = \smfrac{1}{2} \, \alpha / \|X_E\|_E $,   then $ E(t) = \smfrac{1}{1+t} \big( E + t X_E \big) $  lies  in $ \swath( \alpha) \; . $
\end{prop}
\vspace{2mm}

Proposition \ref{prop.bc}  is proven in \S\ref{sect.d}, where it is observed to immediately \vspace{1mm} follow from   \vspace{1mm} Corollary~\ref{cor.db}. 

Proposition \ref{prop.bc}  suggests a simple algorithm: Given an initial matrix $ E_0 $ contained in  $ \swath( \alpha ) $, recursively compute a sequence of matrices in $ \swath( \alpha ) $ with monotonically decreasing objective values, according to
\[    E_{i+1} = \smfrac{1}{1+t_i} \big( E_i + t_i X_i) \textrm{ where $ X_i $ is optimal for $ \qp_{E_i}( \alpha ) $ and $ t_i = \smfrac{1}{2} \, \alpha/ \|X_i\|_{E_i} $} \; . \]
If $ X_i $  happens to satisfy, say, $ \| X_i \|_{E_i} \leq \sqrt{n} $, then according to (\ref{eqn.bc}), 
\[   \frac{\tr(C E_{i+1}) - \textrm{opt\_val} }{\tr(C  E_i) - \textrm{opt\_val}  }  \leq   \frac{1}{1 + \smfrac{1}{2} \, \alpha/ \sqrt{n}  }  \; , \]
which is exactly the kind of inequality that, if it held for all $ i $ and $ n $, and held regardless of the initial matrix $ E_0 \in \swath( \alpha ) $, would yield a complexity bound matching the best one known for semidefinite programming. \vspace{1mm}

But, alas, in general a bound such as $ \| X_i \|_{E_i} \leq \sqrt{n} $ need not hold.  Thus, if we are to succeed in matching the best complexity bound known for semidefinite programming, we should find a relevant way to measure the algorithm's progress even for iterations in which $ \| X_i \|_{E_i} $ is large relative to $ \sqrt{n} $.  And if there is to be any hope of having an algorithm that is interesting in practice as well as in theory, we need to replace the rigid choice $ t = \smfrac{1}{2} \, \alpha / \|X_E\|_E $ with an expression for $ t $ that can potentially evaluate to a far larger value.  \vspace{1mm}

The key to making further progress lies in duality theory.

\section{\bf  Duality, and the Main Theorem for Semidefinite Programming}  \label{sect.c}

Recall that for a convex cone $ K $, the dual cone with respect to an inner product $ \langle \; , \; \rangle $ is defined as
\[  K^* := \{ s: \lin  x, s \rin \geq 0 \textrm{ for all $ x \in K $} \} \; . \]
For $ E \in \Sympp \; , $ let $ K_E(\alpha)^* $ be the cone dual to $ K_E( \alpha) $ with respect to the trace inner product on $ \Sym $.
 \begin{lemma}  \label{lem.ca}
 $ K_E( \alpha )^* = K_{E^{-1}}\big( \sqrt{n - \alpha^2} \, \big) $ \end{lemma}
\noindent {\bf Proof:}  First consider the inner product $ \langle \; , \; \rangle_E \; , $ for which we noted following (\ref{eqn.ba}) 
that the cone $ K_E( \alpha ) $ is circular, consisting of all ``vectors'' $ X \in \Sym $ for which the angle with $ E $ does not exceed $ \arccos( \alpha/ \sqrt{n}) $.  Clearly, for this inner product the dual cone -- which we denote $ K_E( \alpha )^{*_E} $ -- is the circular cone consisting of ``vectors'' $ \bar{S}   \in \Sym $ for which the angle with $ E $ does not exceed
\[  \pi/2  - \arccos( \alpha/\sqrt{n} ) \, = \, \arccos( \sqrt{1 - \alpha^2/n} \, ) \; . \]
Thus, $ K_E( \alpha )^{*_E} = K_E( \sqrt{n - \alpha^2} \, ) \; . $ 
  However, from the definition of dual cones, it is easily checked that 
  \[  \bar{S}  \in K_E( \alpha)^{*_E} \quad \Leftrightarrow \quad E^{-1} \bar{S}  E^{-1} \in K_E( \alpha )^* \; . \]
  Hence, 
\begin{align*}
   K_E( \alpha )^* & = \{ E^{-1} \bar{S} E^{-1}: \bar{S} \in K_E( \alpha )^{*_E} \}  \\  
    & = \{ E^{-1} \bar{S} E^{-1}: \bar{S} \in K_E( \sqrt{n - \alpha^2} \, ) \} \\
    & = \{ E^{-1} \bar{S} E^{-1}:  \tr( E^{-1} \bar{S}) \geq \sqrt{n - \alpha^2} \, \, \tr( \, (E^{-1} \bar{S})^2 \,  )^{1/2} \}   \\
   & = \{ S: \tr(E S) \geq  \sqrt{n - \alpha^2} \, \, \tr( \, (E S)^2 \, )^{1/2}\}   \\
   & = K_{E^{-1}}( \sqrt{n - \alpha^2} \, ) \; ,          
    \end{align*}
    completing the proof.  \hfill $ \Box $
 \vspace{3mm}

The optimization problem dual to $ \qp_E( \alpha) $ is
\[ \left. \begin{array}{rl}
\min_{y,S} & b^T y \\
\textrm{s.t.} & {\mathcal A}^* y + S = C \\
  & X \in K_E( \alpha )^*   \end{array} \right\} \textrm{QP}_E( \alpha)^* \; ,   \]    
recalling that $ {\mathcal A}^* y = \sum_i y_i A_i \,  $.  The optimal solution to this problem plays an extremely important role for us, not only in the analysis, but even in specifying the algorithm in such a way that the rigid choice $ t = \alpha / \|X_E\|_E $ is replaced with an expression for $ t $ that can potentially evaluate to a far larger value. \vspace{1mm}
 
Assume $ E \in \swath( \alpha ) $ and let $ X_E = X_E( \alpha ) $ be the optimal solution for $ \qp_E( \alpha ) $.  We noted in (\ref{eqn.bb}) that among the first-order conditions satisfied by $ X_E $ is there exist $ y \in \mathbb{R}^{m} $ and $ \lambda \in \mathbb{R}  $ for  which
\begin{equation}  \label{eqn.ca}
   0 = \lambda C + \sum_i y_i A_i  +  \tr(E^{-1}X_E) \, E^{-1}  - \alpha^2 E^{-1} X_E E^{-1} \; . \end{equation} 
Standard arguments based on this equation show how to determine the optimal solution of the dual problem $ \qp_E( \alpha )^* $.  We now give the arguments for the sake of readers who might be unfamiliar with the reasoning.  \vspace{1mm}

Multiplying both sides of (\ref{eqn.ca})  on the right by $ E - X_E $, applying the trace operator and using $ \tr\big( A_i \,  (E - X_E) \big) = 0 $ (because $ {\mathcal A} (E) = b = {\mathcal A} (X_E) $), gives upon straightforward simplification that
\begin{align*}
  0 & = \lambda \, \tr\big( C (E - X_E) \big)  + (n - \alpha^2) \, \tr(E^{-1} X_E) - \tr(E^{-1}X_E)^2 + \alpha^2 \tr\big( \, ( E^{-1} X_E)^2 \, \big)    \\
     & = \lambda \, \tr\big( C (E - X_E) \big)  + (n - \alpha^2) \, \tr(E^{-1} X_E) \qquad \textrm{(because $ X_E \in \partial K_E( \alpha ) $)} \; . 
\end{align*}
Thus, 
\begin{equation}  \label{eqn.caa}
  \lambda = - \frac{(n - \alpha^2) \, \tr(E^{-1}X_E) }{ \tr\big( C(E-X_E) \big) } \; . 
  \end{equation} 
Hence, defining 
\begin{equation}  \label{eqn.cb}
  y_E = y_E( \alpha ) := - \smfrac{1}{\lambda } y \quad \textrm{(where $ y $ is as in (\ref{eqn.ca}))} 
  \end{equation} 
and 
\begin{align}
   S_E = S_E( \alpha ) & :=   - \smfrac{1}{ \lambda } \big(  \tr(E^{-1}X_E) \, E^{-1}  - \alpha^2 E^{-1} X_E E^{-1} \big) \nonumber \\
& =  \smfrac{\tr(C  (E - X_E)) }{n - \alpha^2} \big( E^{-1}  - \smfrac{\alpha^2}{ \tr( E^{-1} X_E)}  E^{-1} X_E E^{-1} \big) \; ,  \label{eqn.cc}
\end{align}
we have from (\ref{eqn.ca})  that $ (y_E,S_E)  $ satisfies the linear equations $ {\mathcal A}^* y + S = C $.  To establish feasibility for $ \qp_E( \alpha )^* $, it remains to show $ S_E \in K_E( \alpha )^* $, that is, it remains to show $ S_E \in K_{E^{-1}}\big( \sqrt{n - \alpha^2} \big) $ (making use of Lemma~\ref{lem.ca}).   \vspace{1mm}

However, a simple computation using (\ref{eqn.cc})   shows
\begin{equation}  \label{eqn.cd}
    \tr( E S_E) = \tr\big( C( E - X_E) \big)    \; . \end{equation} 
Likewise, using also $ \tr( \, ( X_E E^{-1})^2 \, ) = \smfrac{1}{ \alpha^2} \tr( E^{-1} X_E) $ (because $ X_E \in \partial K_E( \alpha ) $), 
\begin{equation}  \label{eqn.ce}
  \tr\big( \, (E S_E)^2 \, \big) =  \frac{\tr\big( C( E - X_E) \big)^2   }{n - \alpha^2} \; . \end{equation}
Thus, $ S_E \in K_{E^{-1}} \left(  \sqrt{n - \alpha^2} \,  \right) = K_E( \alpha )^* $.  Hence, $ (y_E, S_E) $ is feasible \vspace{1mm} for $ \qp_E( \alpha )^* $. \vspace{1mm}

In fact, the pair $ (y_E, S_E) $ is optimal for $ K_E (\alpha )^* $.  Indeed, since
 \begin{align*}
\tr( X_E S_E) & =  \smfrac{\tr(C  (E - X_E)) }{n - \alpha^2} \, \tr \left(   X_E \left(  E^{-1}  - \smfrac{\alpha^2}{ \tr( E^{-1} X_E)}  E^{-1} X_E E^{-1} \right) \right) \\
& \smfrac{\tr(C  (E - X_E)) }{(n - \alpha^2) \, \tr( E^{-1} X_E)} \left( \tr(E^{-1} X_E)^2 - \alpha^2 \, \tr( \, ( E^{-1} X_E)^2 \, ) \right)  \\
& = 0 \qquad  \textrm{(because $ X_E \in \partial K_E( \alpha ) $)} \; ,       
\end{align*}
     we have
  \begin{align} 
  b^T y_E & = {\mathcal A} (X_E)^T y_E \nonumber \\
           & =  \tr\big( X_E \, ({\mathcal A}^* y_E)\big)  \nonumber  \\
           & =  \tr\big(  X_E (C - S_E ) \big) \nonumber \\
           & = \tr( C X_E) \; , \label{eqn.cf}
           \end{align}
 showing that the feasible points $ X_E $ and $ (y_E, S_E) $ satisfy strong duality, and hence are optimal for the primal-dual pair of optimization problems $ \qp_E( \alpha ) $, $ \qp_E( \alpha )^* $.
           \vspace{1mm}
           
It is not, however, the pairing of $ X_E $ and $ (y_E, S_E) $ that is of greatest pertinence, but rather, the pairing of $ E $ and $ (y_E, S_E) $.  The reason is that -- assuming $ 0 < \alpha \leq 1 $ -- this is a pairing of feasible points for the primal-dual problems of real interest, $ \sdp $ and $ \sdp^* $.  Indeed,  $ E $ is trivially (strictly) feasible for $ \sdp $, because $ E \in \swath( \alpha ) $.  On the other hand, when $ 0 < \alpha \leq 1 $, we know $ \Symp \subseteq K_E( \alpha ) $ (Lemma~\ref{lem.baa}), and hence, $ K_E( \alpha )^* \subseteq (\Symp)^* = \Symp $.   Since $ {\mathcal A}^* y_E + S_E = C $, it is indeed the case that $ (y_E, S_E) $ is feasible for $ \sdp^* $. \vspace{1mm}

With a primal-dual feasible pair for $ \sdp $ in hand, it now makes sense to refer a duality gap: 
\begin{align}
    \textrm{gap}_E & := \tr( C E ) - b^T y_E \nonumber \\ 
             & = \tr\big( C (E - X_E) \big)    \quad \textrm{(by (\ref{eqn.cf}))} \; .  \label{eqn.cg}
\end{align} 
(Thus, the difference in the primal objective values for the points $ E $ and $ X_E $ is a duality gap in disguise.)  \vspace{1mm}

 In order to match the best complexity bound for semidefinite programming, we want our algorithm to halve the duality gap in $ {\mathcal C}    \sqrt{n} $ iterations for some constant $ {\mathcal C}   $.  The following theorem implies our algorithm accomplishes this with a constant dependent only on $ \alpha $.     \vspace{1mm}

\hypertarget{targ_main_thm_sdp}{}
\begin{main_thm_sdp}   Assume $ 0 < \alpha < 1 $ and $ E \in \swath( \alpha ) $.  Let $ X_E( \alpha) $, $  (y_E( \alpha), S_E( \alpha)) $ be the optimal solutions of $ \qp_E(\alpha) $, $ \qp_E(\alpha)^* $, respectively. Define
\begin{equation}  \label{eqn.caaa}
   E' = \smfrac{1}{1 + t_{E}( \alpha)} \, \left( E + t_{E}( \alpha ) \, X_{E}( \alpha) \right) \; , 
   \end{equation} 
where $  t_{E}(\alpha)  $ minimizes the convex quadratic polynomial 
\[  t \mapsto \tr \bigg(   \, \Big(  \big( E + t \, X_{E}( \alpha) \big) \, S_{E}( \alpha) \Big)^2 \, \bigg) \; . \]
Then 
\begin{itemize}
\item $ E' \in \swath( \alpha )   $  
\item  $ \tr( C \, E ) > \tr( C \, E' )    $ \quad (primal objective monotonicity)
\item $ b^T y_{E}(\alpha) \leq  b^T y_{E'}( \alpha)  $ \quad (dual objective monotonicity) \; . 
\end{itemize}
Moreover, if beginning with $ E_0 \in \swath( \alpha) $, the identity (\ref{eqn.caaa}) is recursively applied to create a sequence $ E_0, E_1, \ldots $, then for every $ i = 0, 1, \ldots \; , $  
\begin{equation}  \label{eqn.ch}
     \frac{ \mathrm{gap}_{E_{j+1}}}{ \mathrm{gap}_{E_{j}}} \leq  1 - \frac{\kappa}{\kappa + \sqrt{n}}     \quad \textrm{for $ j = i $ or $ j = i+1 $ (possibly both)} \; ,   \end{equation} 
 where
 \[ \kappa := \alpha \sqrt{\smfrac{1- \alpha }{8} } \; . \]           
\end{main_thm_sdp}
\vspace{2mm}

The particular choice $ t = t_{E}( \alpha ) $ -- where $ t_{E}(\alpha) $ minimizes the quadratic polynomial -- is not critical.  As will become clear, the proofs apply to any value $ t $ chosen from the (non-empty) interval $ ( t_E( \alpha ) - \delta, t_E( \alpha ) + \delta ) $ where $ \delta :=  t_E( \alpha ) - \smfrac{1}{2} \alpha/\| X_E( \alpha) \|_E \; . $ In some cases, choosing $ t = t_E( \alpha) $ results in a long step.  (Keep in mind that the step $ E(t) - E = \smfrac{t}{1+t} ( X_E( \alpha) - E) $ does not scale linearly in $ t $.)   \vspace{1mm}

For the purpose of matching the best complexity bound known for SDP, it suffices to fix $ \alpha $ independent of $ n $.  Still, the inequality (\ref{eqn.ch}) is stronger if $ \alpha $ is not close to either $ 0 $ or $ 1 $.  This raises the question:  If $ E_0 \in \swath( \alpha_0) $ for $ \alpha_0 $ close to $ 0 $ or $ 1 $, can the algorithm be modified so as to quickly obtain a point in $ \swath( \alpha) $ for user-chosen $ \alpha$ away from 0 and 1? \vspace{1mm}

   Now, if the user-chosen value $ \alpha $ satisfies $ \alpha  \geq  \alpha_0 $, then $ E_0 $ lies not only in $ \swath( \alpha_0 ) $ but also in $ \swath( \alpha  ) $, so no modification of the algorithm is necessary.  (However, applying the computations with $ \alpha  $ rather than $ \alpha_0 $ will result in iterates $ E_i $ guaranteed to lie in $ \swath( \alpha  ) $ but not necessarily in $ \swath( \alpha_0 ) $.) \vspace{1mm}

On the other hand, if $ \alpha < \alpha_0 $, then obtaining a point in $ \swath( \alpha  ) $ is not as trivial.  Nonetheless, we show (Corollary~\ref{cor.dc}) that a slight modification to the procedure gives an algorithm which computes a point in $ \swath( \alpha  ) $ within $ O \left(   \log( \frac{\alpha_0}{ \alpha} ) + \log( \frac{1 - \alpha}{1 - \alpha_0} ) \right)    $ iterations. \vspace{1mm}

The next two sections lay groundwork for proving the main theorem above.  The proof is completed in \S\ref{sect.f}, after which attention switches from semidefinite programming to hyperbolic programming. 
\vspace{1mm}

\section{{\bf  Staying Within the Swath}} \label{sect.d}

This section is mostly devoted to proving the following proposition, which the subsequent corollary  shows to imply that if the algorithm's initial iterate $ E_0 $ lies in $ \swath( \alpha ) $, then the ensuing iterates $ E_1, E_2, \ldots $ lie not only within $ \swath( \alpha) $, but  lie within $ \swath( \beta) $ for a specified value $ \beta < \alpha $.
 \vspace{1mm}

We use the notation $ \int( {\mathcal S}) $ to denote the interior of a set $ {\mathcal S} $. \vspace{1mm}

\begin{prop}   \label{prop.da}
Assume $ E \in \Sympp \; , $ $ 0 < \alpha < 1 $ and $ 0 \neq X  \in \partial K_E( \alpha ) $.  Let  \[   S  := \smfrac{1}{n - \alpha^2}  \big( E^{-1} - \smfrac{\alpha^2}{\tr(E^{-1} X)  }  E^{-1} X  E^{-1}  \big)   \; . \]
The quadratic polynomial
\begin{equation}  \label{eqn.da}
   q(t) :=  \tr \left(   \, \left(  ( E + t X ) \, S  \right)^2 \, \right)  
   \end{equation} 
is strictly convex and its  minimizer $ \bar{t}  $ satisfies $ \delta :=   \bar{t} -\smfrac{1}{2} \alpha / \| X  \|_E > 0 $.  Moreover,
\[  \bar{t} - \delta \leq t \leq \bar{t} + \delta \quad \Rightarrow \quad E(t) \in \Sympp \, \textrm{ and } \, S \in \int\big( K_{E(t)}( \beta )^* \big) \]
where $ E(t) = \smfrac{1}{1+t} (E + tX) $ and $  \beta = \alpha \sqrt{\smfrac{1 + \alpha }{2} } \; . $ 
\end{prop} 
\vspace{2mm}

Recall that for $ E \in \swath( \alpha ) $, we use the notation $ X_E( \alpha ) $ to denote the optimal solution for the quadratic-cone optimization problem $ \qp_E( \alpha) $, and $ (y_E( \alpha), S_E( \alpha )) $ to denote the optimal solution for the dual problem, $ \qp_E( \alpha)^* \; . $   Following is the consequence of Proposition~\ref{prop.da} that will be used in proving the Main Theorem for Semidefinite Programming.   \vspace{3mm}

\begin{cor}  \label{cor.db}
Assume $ 0 < \alpha < 1 $ and $ E \in \swath( \alpha ) \; . $ Let $ X_E = X_E( \alpha ) $ and $ S_E = S_E( \alpha ) $.  The quadratic polynomial
\begin{equation}  \label{eqn.db}
     t \mapsto \tr \left( \, \big( ( E + t X_E) \, S_E \big)^2 \, \right) 
     \end{equation} 
is strictly convex and its minimizer $ t_E = t_E( \alpha) $ satisfies $ \delta := t_E - \smfrac{1}{2} \alpha/ \|X_E \|_E > 0 $.  Moreover,
\[  t_E - \delta \leq t \leq t_E + \delta \quad \Rightarrow \quad E(t) \in \swath( \beta ) \, \textrm{ and } \, S_E \in \int \big( K_{E(t)}( \beta )^* \big) \; ,  \]
where $ E(t) := \smfrac{1}{1+t} (E + t X_E) $ and $  \beta = \alpha \sqrt{\smfrac{1 + \alpha }{2} }  \; . $ 
\end{cor}  
\noindent {\bf Proof:}  Choosing $ X = X_E $, then $ S $ in Proposition~\ref{prop.da} is a  positive scalar multiple of $ S_E $ (by (\ref{eqn.cc})), and the polynomial (\ref{eqn.db}) is a positive multiple of the polynomial (\ref{eqn.da}).  Consequently, Proposition~\ref{prop.da} shows
\[  t_E - \delta \leq t \leq t_E + \delta \quad \Rightarrow \quad E(t) \in \Sympp \, \textrm{ and } \, S_E \in \int \big( K_{ E(t)}( \beta )^* \big) \; , \]
where $ \beta = \alpha \sqrt{\frac{1 + \alpha}{2} } \; . $ In particular, the quadratic-cone optimization problem $ \qp_{E(t)}( \beta ) $ and its dual problem both are strictly feasible (indeed, $ E(t) $ is strictly feasible for $ \qp_{E(t)}( \beta ) $ and $ (y_E, S_E) = (y_E( \alpha), S_E( \alpha )) $ is strictly feasible for the dual problem).  Consequently, $ \qp_{E(t)}(\beta ) $ has an optimal solution, i.e., $ E(t) \in \swath( \beta ) \; . $ \hfill $ \Box $
 \vspace{3mm}
 
In initially motivating our algorithm in section \ref{sect.b}, we relied upon Proposition \ref{prop.bc}, which asserts that if $ E \in \swath(\alpha) $ and $ t = \smfrac{1}{2} \, \alpha/\|X_E( \alpha)\|_E $, then $ E(t) = \smfrac{1}{1+t} (E + t X_E) $ lies in $ \swath( \alpha) $, where $ X_E = X_E( \alpha) $.  We are now in position to establish the validity of the proposition.   Indeed, it trivially follows from Corollary~\ref{cor.db} by choosing $ t = t_E - \delta = \smfrac{1}{2} \, \alpha/ \|X_E\|_E $. \vspace{1mm}

 Before proceeding to prove Proposition~\ref{prop.da}, we record a corollary that establishes a claim made just after the statement of the Main Theorem for Semidefinite Programming.  Specifically, we claimed  
 a slight modification in the computational procedure that had been described would result in an algorithm which given a point $ E_0 \in \swath( \alpha_0) $, computes a point in $ \swath( \alpha ) $ (where $ \alpha < \alpha_0 $ is user-chosen)  in $  O \left(   \log( \smfrac{\alpha_0}{\alpha}) + \log( \frac{1 - \alpha}{ 1 - \alpha_0}) \right)  $    iterations.   The claim is now easily justified. 

\begin{cor}  \label{cor.dc}
 Assume $ 0 < \alpha < \alpha_0 < 1 $ and $ E_0 \in \swath( \alpha_0 ) $.   Recursively define
\begin{equation}  \label{eqn.dc}
        \alpha_{i+1} = \alpha_i \sqrt{\smfrac{1 + \alpha_i}{2} } \quad \textrm{and} \quad  E_{i+1} = \smfrac{1}{1 + t_i } 
   ( E_i + t_i X_i ) \; ,
        \end{equation} 
 where
 \[  X_i = X_{E_i}( \alpha_i) \quad \textrm{and} \quad t_i = \smfrac{1}{2} \alpha_i/ \|X_i\|_{E_i} \; . \]
Then
\[  i \geq \smfrac{2}{\log(8/7)} \,  \log( \smfrac{\alpha_0}{ \alpha })  + \smfrac{1}{\log(9/8)} \,  \log( \smfrac{1- \alpha }{ 1 - \alpha_0 })  \quad \Rightarrow \quad  E_i \in \swath( \alpha  ) \; . \]     
\end{cor}
\noindent {\bf Proof:} 
Induction using Corollary~\ref{cor.db}  shows $ E_i \in \swath( \alpha_i) $ for all $ i = 0, 1, \ldots  $  .  Thus, our goal is merely to show 
\[ i \geq  \smfrac{2}{\log(8/7)} \,  \log( \smfrac{\alpha_0}{ \alpha })  + \smfrac{1}{\log(9/8)} \,  \log( \smfrac{1- \alpha }{ 1 - \alpha_0 })  \quad \Rightarrow \quad  \alpha_i \leq \alpha  \; . \]   
For this it suffices to show there exists a value $ 0 < \bar{\alpha} < 1 $  such that both
\begin{equation}  \label{eqn.dd}
  \alpha_i \leq \bar{\alpha}  \quad \Rightarrow \quad \frac{ \alpha_{i}}{ \alpha_{i+1}} \geq \sqrt{8/7 } 
  \end{equation}
  and
\[ 
   \alpha_i \geq \bar{\alpha}  \quad \Rightarrow \quad \frac{1 - \alpha_{i+1}}{1 - \alpha_i} \geq 9/8  \; . 
\] 
We show these implications hold for $ \bar{\alpha} = 3/4 $.  
\vspace{1mm}

Indeed, when $ \bar{\alpha} = 3/4 $, the implication (\ref{eqn.dd}) is trivial from the identity for $ \alpha_{i+1} $ in    (\ref{eqn.dc}).   
   On the other hand, squaring both sides of that identity, then substituting $ \alpha_i = 1 - \epsilon_i $ and $ \alpha_{i+1} = 1 - \epsilon_{i+1} $, gives   
\begin{align*}
   2 \epsilon_{i+1} - \epsilon_{i+1}^2 & = \smfrac{5}{2} \epsilon_i - 2 \epsilon_i^2 + \smfrac{1}{2} \epsilon_i^3 \\
& >  \smfrac{5}{2} \epsilon_i - 2 \epsilon_i^2 \; . 
\end{align*}
However, $ \epsilon_i < \epsilon_{i+1} $, and hence,
\[      2 \epsilon_{i+1} > \smfrac{5}{2} \epsilon_i - \epsilon_i^2 \; . \]
Thus, 
\[    \epsilon_i \leq \smfrac{1}{4} \quad \Rightarrow  \quad \epsilon_{i+1} \geq  \smfrac{1}{2} \, ( \smfrac{5}{2} - \smfrac{1}{4} ) \, \epsilon_i = \smfrac{9}{8} \, \epsilon_i \; , \]
that is,
\[ \alpha_i \geq 3/4 \quad \Rightarrow \quad \frac{1 - \alpha_{i+1}}{1 - \alpha_i} \geq 9/8 \; , \]
completing the proof. \hfill $ \Box $
   \vspace{3mm}

Now we turn to the proof of Proposition \ref{prop.da}, which relies on a sequence of results.  For ease of reference, we again display the notation and the assumptions appearing in Proposition~\ref{prop.da}:
\begin{gather} 
    E \in \Sympp \; , \quad 0 < \alpha < 1 \; , \quad 0 \neq X \in \partial K_E( \alpha) \; , \label{eqn.df}\\
  \quad S := \smfrac{1}{(n - \alpha^2)} (E^{-1} - \smfrac{\alpha^2}{\tr( E^{-1} X) } E^{-1} X E^{-1})  \; ,  \label{eqn.dg} \\
   q(t) :=  \tr \Big( \, \big( (E + tX) \, S \big)^2 \, \Big)  \; .   \label{eqn.dh}
   \end{gather}  
We begin the proof of Proposition~\ref{prop.da} with a simple lemma. \vspace{1mm}

\begin{lemma}  \label{lem.dd}  If $ E $, $ X $ and $ S $ are as (\ref{eqn.df}) and (\ref{eqn.dg}), then 
\[  \tr( E S) = 1 \quad \textrm{and} \quad   \tr(XS) = 0 \; . \] 
\end{lemma}
\noindent {\bf Proof:}  Simply observe that
\begin{align*}
\tr( E S) & = \smfrac{1}{n - \alpha^2} \, \tr\big( E ( E^{-1} - \smfrac{\alpha^2}{\tr(E^{-1} X) } E^{-1} X E^{-1}) \big) \\ 
& =  \smfrac{1}{n - \alpha^2} \, \big( n -  \smfrac{\alpha^2}{\tr(E^{-1} X) } \tr( X E^{-1} ) \big) \\
& = 1 \; , 
\end{align*}
and      
\begin{align*}
\tr(XS) & = \smfrac{1}{n - \alpha^2} \, \tr \big(  X ( E^{-1} - \smfrac{\alpha^2}{\tr(E^{-1} X) } E^{-1} X E^{-1}) \big) \\
& = \smfrac{1}{(n - \alpha^2) \, \tr(E^{-1}X) } \Big( \tr(E^{-1} X)^2 - \alpha^2 \, \tr \big( \, ( E^{-1} X)^2 \, \big) \Big) \\
& = 0 \qquad \textrm{(because $ X \in \partial K_E( \alpha) $)} \; , 
 \end{align*}     
thus completing the proof. \hfill $ \Box $
\vspace{3mm}

The next proposition gives meaning to the polynomial $ q$ defined by (\ref{eqn.dh}).  Two elements in this proposition's proof are completed only in the proof of the subsequent proposition, i.e., Proposition~\ref{prop.df}.  The proof of Proposition~\ref{prop.df} is self-contained, and so it might seem the natural order would be to present Proposition~\ref{prop.df} first rather than second.  Our reason for switching the order is to provide motivation for the polynomial $ q $ before displaying the string of technical calculations needed for the proof of Proposition~\ref{prop.df}. \vspace{1mm}

\begin{prop}   \label{prop.de}
Let $ E $, $ \alpha $, $ X $ and $ S $ and $ q $ be as in (\ref{eqn.df}), (\ref{eqn.dg}) and (\ref{eqn.dh}).
For every $ 0 < \beta < 1 $ and $ t > -1 $,  
\[ \left(  E(t) \in \Sympp \right)  \,  \wedge \,   \left(  S \in \int(K_{E(t)}( \beta )^*) \right) \quad \Leftrightarrow \quad   q(t) <  \frac{1}{n - \beta^2}  \]
(where $ E(t) = \smfrac{1}{1+t} (E + tX) \, $). 
Furthermore, $ q $ is strictly convex, $ q(0) = 1/(n - \alpha^2) $ and $ q'(0) < 0 $.
\end{prop}
\noindent {\bf Proof:}  For any fixed $ \bar{S} \in \Sympp \; , $ begin by considering the function 
\begin{equation}  \label{eqn.di}
   \bar{E}  \mapsto \frac{ \tr \left( \, \left( \bar{E}  \bar{S}  \right)^2 \right)^{1/2}}{ \tr( \bar{E}  \bar{S} ) } \; . 
   \end{equation} 
The function is well-defined if we restrict  $ \bar{E}  $ to lie in $  \Symp \setminus \{ 0 \}  $ -- indeed, the value of the denominator is then $ \sum_i \lambda_i > 0 $, where $ \lambda_1, \ldots, \lambda_n $ are the eigenvalues of $ \bar{S}^{-1/2} \bar{E}  \bar{S}^{-1/2} \in  \Symp \setminus \{ 0 \} \; . $  \vspace{1mm}

The value of the numerator is $ \left(   \sum_i \lambda_i^2 \right)^{1/2} $, and thus by the relation between the 1-norm and the 2-norm, the value of the function for $ \bar{E}  \in \Symp \setminus \{ 0 \}  $ satisfies
\[   \frac{ \left( \sum_i \lambda_i^2 \right)^{1/2}}{ \sum_i \lambda_i } \geq \frac{1}{\sqrt{m} } \; , \]
where $ m $ is the number of nonzero eigenvalues.  In particular, the value of the function is at least $ 1/ \sqrt{n}  $ for all $ \bar{E}  \neq 0 $ in the set $  \Symp  $, and has value at least $ 1/ \sqrt{n-1}  $ for all $ \bar{E}  \neq 0 $ in the boundary. \vspace{1mm}

Now observe that for $ \bar{E} \in \Sympp, $ the value of the function (\ref{eqn.di})  is equal to $ 1/ \gamma $ for the largest value $ \gamma $ satisfying $ \bar{S}  \in K_{ \bar{E}^{-1}}( \gamma ) $ -- equivalently, is equal to $ 1/\sqrt{  n - \beta^2} $ for the smallest non-negative $ \beta  $ satisfying $ \bar{S}  \in K_{ \bar{E}^{-1}}(   \sqrt{  n - \beta^2} \, ) = K_{ \bar{E}}  (\beta)^* $ (Lemma~\ref{lem.ca}). \vspace{1mm}

Substituting $ \bar{E} = \smfrac{1}{1+t} (E + tX) $ in (\ref{eqn.da}), and also $ \bar{S} = S $, gives   
\begin{align*}
& \frac{ \tr \left( \, \left( \smfrac{1}{1+t} ( E + t X) \, S \right)^2 \right)^{1/2}}{ \tr( \, \smfrac{1}{1+t} (E + t X) \, S ) } \\ & = \frac{ \tr \left( \, \left(  ( E + t X) \, S \right)^2 \right)^{1/2}}{ \tr( \,  (E + t X) \, S ) } \\
& =  \tr \left( \, \left(  ( E + t X) \, S \right)^2 \right)^{1/2} \qquad  \textrm{(by Lemma~\ref{lem.dd})}  \\
& = \sqrt{q(t)} \; .   
\end{align*}
Thus, we know that if $ E(t) \in \Sympp $, then $ q(t) = 1/( n - \beta^2) \; , $ where $ \beta $ is the smallest non-negative value for which $ S $ is contained in $  K_{E(t)}( \beta )^* \; . $ Moreover, we know that if $ E(t) $ lies in the boundary of $ \Symp $, then $ q(t) \geq 1/(n-1) \; . $ \vspace{1mm}

To complete the proof, it now suffices to show: 
\begin{quote} 
 {\bf A:} \,  For all $ 0 < \beta < 1 $ and $ t > - 1 \; , $  \\
$ \textrm{~} $ \qquad  $  q(t) \leq 1/(n - \beta^2) \, \, \Rightarrow \, \,   E(t) \in \Sympp  $ \\
 {\bf B:} \, The quadratic polynomial $ q $ is strictly convex \\
{\bf C:} \, $ q(0) = 1/(n - \alpha^2) $  \\
{\bf D:} \,  $ q'(0) < 0 $ 
 \end{quote}
 Of course property {\bf A}  is equivalent to the property that if  $ E(t) \notin \Sympp $ and $ t > -1 \; , $ then $ q(t) \geq 1/( n-1) \; . $ However, we already know the inequality $ q(t) \geq 1/(n-1) $  to hold if $ E(t) $ lies in the boundary of $ \Symp $. Thus, since $ E(0) = E \in \Sympp \; , $ to establish property {\bf A}  it suffices to establish property {\bf B} and $ q(0) < 1/(n-1) $; hence, it suffices to establish properties {\bf B} and {\bf C}.   \vspace{1mm}

 In turn, since $ q(t) \geq 1/(n-1) $ for some $ 0 < t < 1 $ (in particular, for the choice of $ t $ for which $ E(t) $ is in the boundary of $ \Symp $), to establish property {\bf B} for the \underline{quadratic} polynomial $ q $, it suffices to establish properties {\bf C} and  {\bf D}. We defer establishing properties {\bf C} and {\bf D}  until the the proof of Proposition~\ref{prop.df}, in order to avoid duplicating calculations.  \hfill $ \Box $ \vspace{3mm}

\begin{prop}  \label{prop.df}
If $ E $, $ \alpha $, $ X $ and $ S $ and $ q $ are as in (\ref{eqn.df}), (\ref{eqn.dg}) and (\ref{eqn.dh}), then 
\[ 0 < t \leq \alpha / \|X \|_E \quad \Rightarrow \quad q(t) < \smfrac{1}{ n - \alpha^2} \,  \left( 1 - 2 \, t \, \smfrac{1 - \alpha }{n - \alpha^2}  \,  \|X\|_E \,  \left( \alpha - t \, \|X \|_E \right)   \,\right)  \; .  \]
\end{prop}
\noindent {\bf Proof:}  Let $ \bar{E}(t) :=  E + t X $ ($ = (1+t) E(t) $), and to avoid having fractions in most calculations, define
\[    \bar{S} = (n - \alpha^2) \, \tr(E^{-1} X) \,  S =  \tr(E^{-1} X) \, E^{-1} - \alpha^2 E^{-1} X E^{-1} \; . \]
Thus,
\begin{equation}     \label{eqn.dj}
   q(t) =    \frac{1}{ (n - \alpha^2)^2 \, \tr( E^{-1} X)^2  }   \, \tr \left( \, \left( \bar{E}(t) \, \bar{S} \right)^2 \, \right) \; . \end{equation}
   
To ease notation, let $ Z := X E^{-1} $, whose eigenvalues are the same as for  $ E^{-1} X $ (and also for $ E^{-1/2} X E^{-1/2} $ (hence all eigenvalues are real)).  The  assumption $ t \leq \alpha / \|X\|_E $ is equivalent to $ t \leq \alpha / \tr(Z^2)^{1/2} $. We make extensive use of the fact that $ \tr(Z) = \alpha \,  \tr(Z^2)^{1/2} $, that is, $ \tr( E^{-1} X) = \alpha \, \tr( \, ( E^{-1} X)^2 \, )  $ (i.e., $ X $ is in the boundary of $ K_E( \alpha ) $, by assumption).    \vspace{1mm}
 
We have
\begin{align*}
   \left( \bar{E}(t) \, \bar{S} \right)^2  & = \tr(Z)^2 I \\
   &    \qquad + 2 \, \tr(Z) \, \big( t \, \tr(Z) - \alpha^2 \big) \, Z \\
  & \qquad  \quad + \big( t^2 \, \tr(Z)^2 - 4 t \alpha^2 \, \tr(Z) + \alpha^4 \big) \,  Z^2 \\
  & \qquad  \qquad   - 2t \alpha^2 \, \big( t \, \tr( Z) - \alpha^2 \big) \, Z^3 \\
 & \qquad  \qquad  \quad + t^2 \alpha^4 \,  Z^4 \; .        
   \end{align*}
   Thus,
\begin{align*}
  \tr\big( \, \left( \bar{E}(t) \, \bar{S} \right)^2 \, \big) &  = (n - 2 \alpha^2) \, \tr(Z)^2 + \alpha^4 \, \tr(Z^2) \\
    & \qquad  + 2 t \,    \Big(  \tr(Z)^3 - 2 \alpha^2 \, \tr(Z) \,  \tr(Z^2)  + \alpha^4 \, \tr(Z^3) \Big)      \\
& \qquad  \quad + t^2 \, \Big( \tr(Z)^2 \, \tr(Z^2) - 2 \alpha^2 \, \tr(Z) \, \tr(Z^3)  + \alpha^4 \, \tr(Z^4) \Big)    \; ,   
  \end{align*}
which upon substitution of $ \tr(Z)^2 $ for  $ \alpha^2 \, \tr(Z^2) $ becomes 
\begin{align}
\tr\big( \, \left( \bar{E}(t) \, \bar{S} \right)^2 \, \big) &  = (n - \alpha^2) \, \tr(Z)^2 \label{eqn.dja} \\
    & \qquad  - 2 t \,    \Big( \tr( Z)^3 -  \alpha^4 \, \tr(Z^3)   \Big)   \nonumber   \\
& \qquad  \quad + t^2 \, \Big( \tr(Z)^2 \, \tr(Z^2) - 2 \alpha^2 \, \tr(Z) \, \tr(Z^3)  + \alpha^4 \, \tr(Z^4) \Big)      \; . \nonumber
  \end{align}
(In passing, we note that substituting $ t = 0 $ immediately shows $ q(0) = 1/(n - \alpha^2) $, one of two facts stated in Proposition~\ref{prop.de}  whose proofs were deferred.)  \vspace{1mm}
 
The coefficient of the linear term for $ t $ in (\ref{eqn.dja}) is
\[   - 2  \Big( \tr( Z)^3 -  \alpha^4 \, \tr(Z^3)   \Big)  = - 2 \alpha^3 \Big( \tr(Z^2)^{3/2}  -  \alpha \, \tr(Z^3)   \Big)  \; , \]
which is negative since $ \tr( Z^3) \leq \tr( Z^2)^{3/2} $, $ 0 < \alpha < 1 $ and $ Z \neq 0 \; . $ (Thus, $ q'(0) < 0 $, the remaining fact stated in Proposition~\ref{prop.de}  for which the proof was deferred.)   \vspace{1mm}

On the other hand, from
\[     \tr(Z)^2 \, \tr(Z^2) = \alpha^2 \,  \tr(Z^2)^2 \; , \]
\[        2 \alpha^2 \,  \tr(Z) \, \tr(Z^3) =  2 \alpha^3 \, \tr(Z^2)^{1/2} \, \tr(Z^3)  \]
and
\[     \alpha^4 \tr(Z^4) \leq \alpha^4 \tr( Z^2)^2 \; , \]
we see that the coefficient for the quadratic term for $ t $ in (\ref{eqn.dja}) can be bounded above by
\begin{align*}
 &   \alpha^2 \,  \tr(Z^2)^{1/2} \, \left( (1 + \alpha^2) \, \tr(Z^2)^{3/2} - 2 \alpha \, \tr(Z^3) \right) \\
&  \qquad  <  2 \alpha^2 \, \tr(Z^2)^{1/2} \,\left(  \tr(Z^2)^{3/2} - \alpha \, \tr(Z^3) \right) \; . 
\end{align*}

Thus, substituting for both the linear and quadratic coefficients in  (\ref{eqn.dja}), we have for all $ t \neq 0 $ that
\begin{align*}
\tr\big( \, \left( \bar{E}(t) \, \bar{S} \right)^2 \, \big) & < (n - \alpha^2) \, \tr(Z)^2  \\
    & \qquad  - 2 \, t \, \alpha^3 \,  \Big( \tr(Z^2)^{3/2}  -  \alpha \, \tr(Z^3)   \Big)      \\
& \qquad  \quad + 2 \, t^2 \,   \alpha^2  \, \tr(Z^2)^{1/2} \,\left(  \tr(Z^2)^{3/2} - \alpha \, \tr(Z^3) \right)  \\
& = (n - \alpha^2) \, \tr(Z)^2 - 2 \, t \, \alpha^2 \left( \alpha - t \, \tr(Z^2)^{1/2} \right) \,\left(  \tr(Z^2)^{3/2} - \alpha \, \tr(Z^3) \right) \; .     
  \end{align*}
Hence,  for $ t \neq 0 $
\begin{align*}
q(t) & = \frac{1}{(n - \alpha^2)^2 \, \tr(Z)^2} \, \tr  \big( \, \left( \bar{E}(t) \, \bar{S} \right)^2 \, \big) \qquad  \textrm{(by (\ref{eqn.dj}))} \\  
& < \frac{1}{n - \alpha^2} - \frac{2 t \alpha^2}{ (n - \alpha^2)^2 \, \tr(Z)^2 }  \,   \left( \alpha - t \, \tr(Z^2)^{1/2} \right) \,\left(  \tr(Z^2)^{3/2} - \alpha \, \tr(Z^3) \right) \; .   
\end{align*}
However, upon substituting $ \tr(Z)^2 = \alpha^2 \, \tr(Z^2) $ into the denominator, the last expression is seen to equal  
\[ \frac{1}{n - \alpha^2} - \frac{2 t \, \tr(Z^2)^{1/2} }{ (n - \alpha^2)^2  }  \,   \left( \alpha - t \, \tr(Z^2)^{1/2} \right) \,\left(  1 - \alpha \, \frac{ \tr(Z^3)}{ \tr(Z^2)^{3/2} }  \right) \; . \]
Since $ \tr(Z^3) \leq \tr(Z^2)^{3/2} $, it follows that for all $ 0 < t \leq \alpha/ \tr( Z^2)^{1/2} $ ($ = \alpha / \|X \|_ E $),
\[ q(t) < \smfrac{1}{ n - \alpha^2} \left( 1 - 2 t \,  \smfrac{1 - \alpha }{n - \alpha^2} \, \|X \|_E \,  \left( \alpha - t \, \|X \|_E \right) \right) \; , \]
thus concluding the proof.  \hfill $ \Box $
 \vspace{3mm}

\noindent {\bf Proof of Proposition~\ref{prop.da}:}   
The polynomial $ q $ is strictly convex by Proposition~\ref{prop.de}.   
To verify the minimizer $ \bar{t} $ satisfies $ \bar{t}  > \smfrac{1}{2} \alpha / \| X \|_E $, first apply Proposition~\ref{prop.df} to confirm that 
\[    q( \alpha/ \|X \|_E) < 1/(n - \alpha^2) \; . \]
Since $ q(0) = 1/(n - \alpha^2) $ (Proposition~\ref{prop.de}), the minimizer of the strictly-convex quadratic polynomial $ q $ must thus lie closer to $ \alpha/\|X \|_E $ than to $ 0 $, that is, $ \delta :=  \bar{t} - \smfrac{1}{2} \alpha/\|X \|_E  > 0 $.  \vspace{1mm}

Since $ q $ is a convex quadratic polynomial, its values in any interval centered at its minimizer $ \bar{t} $ do not exceed its value at the endpoints.  Thus,  
\[    \bar{t} - \delta \leq t \leq \bar{t} + \delta \quad \Rightarrow \quad q(t) \leq q(\smfrac{1}{2} \alpha/\|X \|_E) \; . \]
However,  Proposition~\ref{prop.df} easily implies 
\begin{align*}  
 q(\smfrac{1}{2} \alpha/\|X \|_E)  & <    \smfrac{1}{ n - \alpha^2} \,  \left( 1 - \smfrac{1}{2} \,  \alpha^2 \,  \smfrac{1- \alpha }{n - \alpha^2}  \right) \\
& <  \frac{1}{ n - \alpha^2} \, \, \, \,  \frac{1}{1 +  \smfrac{1}{2} \, \alpha^2 \,  \frac{1- \alpha }{n - \alpha^2} } \\
& = \frac{1}{n - \beta^2} \; ,  
 \end{align*} 
where 
\[ \beta = \alpha \sqrt{\smfrac{1 + \alpha }{2}} \; . \]
Thus,
\[ \bar{t} - \delta \leq t \leq \bar{t} + \delta \quad \Rightarrow \quad q(t) < \smfrac{1}{ n - \beta^2} \quad \textrm{for $ \beta = \alpha \sqrt{\smfrac{1+ \alpha }{2}}$} \; . \]
Hence, by Proposition~\ref{prop.de},
\[  \bar{t} - \delta \leq t \leq \bar{t} + \delta \quad \Rightarrow \quad E(t) \in \Sympp  \textrm{ and }    S \in \int(K_{E(t)}( \beta )^*) \; , \]
where $ E(t) = \smfrac{1}{1+t} ( E + tX) $.  \hfill $ \Box $
   \vspace{3mm}
 
 In concluding the section, we record a result that will not play a role until \S\ref{sect.j}, where the focus is hyperbolic programming.  This result is an immediate corollary to calculations done in the proof of Proposition~\ref{prop.df}.    \vspace{1mm}

\begin{cor}  \label{cor.dga}
Assume $ E $, $ \alpha $, $ X $, $ S $ and $ q $ are as in (\ref{eqn.df}), (\ref{eqn.dg}) and (\ref{eqn.dh}).  If $ \lambda_1, \ldots, \lambda_n $ are the eigenvalues of $ E^{-1/2} X E^{-1/2} $,  then the polynomial $ q $ is a positive multiple of 
\[    t \mapsto a t^2 + b t + c \]
where
\begin{gather*}  
a = \Big( \sum_j \lambda_j \Big)^2 \sum_j \lambda_j^2 - 2 \alpha^2 \Big( \sum \lambda_j \Big) \sum_j \lambda_j^3 + \alpha^4 \sum_j \lambda_j^4 \; , \\
 b = 2 \alpha^4 \sum_j \lambda_j^3 - 2 \, \Big( \sum_j \lambda_j \big)^3  \quad \textrm{and} \quad c = (n - \alpha^2) \, \Big( \sum_j \lambda_j \Big)^2 \; .
\end{gather*}    
\end{cor}
\noindent {\bf Proof:}  Immediate from (\ref{eqn.dj}) and (\ref{eqn.dja}). \hfill $ \Box $

\section{{\bf  Conic Geometry}}  \label{sect.e}

Here we develop a result regarding circular cones in Euclidean space.  
Since we rely only on the properties of Euclidean geometry, the results and proofs in this section are written for $ \mathbb{R}^n $ equipped with the dot product and its induced norm, the standard Euclidean norm. Of course the results apply to any specified inner product on a finite-dimensional real vector space, simply because by choosing a basis that is orthonormal for the inner product,  the dot product of the resulting coordinate system is precisely the specified inner product.\vspace{1mm}

For $ e \in \mathbb{R}^n $ satisfying $ \| e \| = \sqrt{n} $ (e.g., the vector of all ones), and for $ 0 < \alpha < \sqrt{n} $, define  
\[     K_e( \alpha ) := \{ x: e^T x \geq \alpha \, \| x \| \} \; , \]
a cone which is invariant under rotations around $ \{ t \, e: t \in \mathbb{R} \} $, and whose dual is 
\begin{align*}
K_e( \alpha )^* & := \{ s: x^T s \geq 0 \textrm{ for all $ x \in K_e( \alpha )$} \} \\
      & = \{ s : e^T s \geq \sqrt{n - \alpha^2}  \, \| s \| \} \; . 
\end{align*}

For subspaces $ L $, we use the notation $ e + L $ to denote the affine space $ \{ e + v : v \in L \} \; . $  \vspace{1mm}

Here is the result to which the section is devoted:
\begin{prop}  \label{prop.ea}
Assume $ e $ satisfies $ \| e \| = \sqrt{n} $, and  $ 0 < \beta \leq \alpha < \sqrt{n} \; .  $ Assume $ \bar{s}  $ lies in the interior of $ K_e( \beta )^* $, and assume $ \bar{x}  $ solves the optimization problem 
\[ 
    \begin{array}{rl}
 \min_x & \bar{s}^T x \\
\textrm{s.t.} & x \in e + L \\
& x \in K_e( \alpha ) \; , \end{array} 
\] 
where $ L $ is a subspace satisfying $ e \notin L $ and $ \bar{s} \notin L^{\perp} $.

Then there exists $ \bar{s}' \in ( \bar{s}  + L^{\perp}) \cap K_e( \alpha )^* $ satisfying 
\[   e^T ( \bar{s}  - \bar{s}' ) \geq {\mathcal C}_1 {\mathcal C}_2   \, \| \bar{s} \| \; ,  \]
where
\[ {\mathcal C}_1 := \left(  \frac{(n - \alpha^2) \| \bar{x} \|^2  \Big( 1 - \big(\frac{\beta}{\alpha} \big)^2 \Big)}{n - \alpha^2 + ( \| \bar{x} \| - \alpha )^2 \Big( 1 - \big(\frac{\beta}{\alpha} \big)^2 \Big)} \right)^{1/2} \]
and
\[ {\mathcal C}_2 :=  \smfrac{1}{n }  \left( \alpha  \, \sqrt{  n - \beta^2} - \beta  \, \sqrt{  n - \alpha^2} \right)   \; . \]
\end{prop}
\vspace{2mm}

The proof of the proposition relies on three lemmas. 
\vspace{1mm}

   \begin{lemma}  \label{lem.eb}
   Assume $ 0 < \beta  \leq \alpha < \sqrt{n}  $, and assume $ e  $ satisfies $ \| e \| = \sqrt{n} $.  Let $ \bar{x} $ be any point in   $  \partial K_e( \alpha  ) $ other than the origin. Assume $ v \neq 0 $ lies in the tangent space to $ \partial K_e( \alpha  ) $ at $ \bar{x} $, as well as lies in the tangent space to $ \partial K_e( \beta  ) $ at some point other than the origin.    Then $ \theta $, the angle between $ \bar{x} $ and $ v $, satisfies $ | \cos ( \theta ) | \leq  \beta  / \alpha  \; . $
\end{lemma}
\noindent {\bf Proof:}  Let $ z \neq 0 $ be the point in $  \partial K_e( \beta  ) $ in whose tangent space $ v $ lies.  Since the line $ \{ z + tv: t \in \mathbb{R} \} $ does not intersect the interior of $ K_e( \beta  ) $, neither $ v $ or $ - v $ is contained in the interior, and thus $ \theta' $, the angle between $ v $ and $ e $, satisfies 
\begin{equation}   \label{eqn.ea}
  | \cos ( \theta' ) | \leq  \beta  / \sqrt{n} 
  \end{equation}
(because $ K_e( \beta ) $ consists of the vectors whose angle with $ e $ does not  \vspace{1mm} exceed $ \arccos (\beta / \sqrt{n}) $).

The tangent space to $ \partial K_e( \alpha  ) $ at $ \bar{x} $ consists of the vectors $ u $ satisfying $  (e^T \bar{x}) (e^T u) - \alpha^2 \, \bar{x}^T u = 0  $, that is, satisfying 
\[ \frac{\bar{x}^T u}{ \| \bar{x} \|} = \frac{e^T u}{ \alpha  } \quad \textrm{(using $ e^T \bar{x} = \alpha  \|\bar{x}\| $)} \; . \]
    Since $ v $ is in this tangent space, $ \theta $ (the angle between $ \bar{x} $ and $ v $) thus satisfies
\[ 
  | \cos ( \theta ) | = \frac{| \bar{x}^T v|}{\|\bar{x}\| \, \|v\|} = \frac{|e^T v|}{ \alpha  \, \| v \|}   = \frac{ \| e \|}{ \alpha  } \, | \cos ( \theta' ) | \leq   \frac{ \beta }{ \alpha  }  \; ,
\] 
where the inequality is due to (\ref{eqn.ea})  and $ \| e \| = \sqrt{n} $.  \hfill $ \Box $
 \vspace{3mm}

\begin{lemma}   \label{lem.ec}
Assume $ 0 < \alpha < \sqrt{n}  $, and assume $ e  $ satisfies $ \| e \| = \sqrt{n} $. Let $ L $ be a subspace, and assume $ \bar{x} \neq 0 $ is some point contained in  $  M := (e + L) \cap \partial K_e( \alpha ) \; . $ Let $ T_{\bar{x}} $ denote the tangent space to $ M $  at $ \bar{x} $, and let $ \theta $ be the angle between $ \bar{x} $ and $ T_{\bar{x}} $ (i.e., the angle between $ \bar{x} $ and its projection onto $ T_{\bar{x}} $).  Then the projection $ P_{L^{\perp}}(e) $ of $ e $ onto $ L^{\perp} $ satisfies
\begin{equation}  \label{eqn.eb}
    \| P_{L^{\perp}}(e)\|^2 = \frac{(n  -  \alpha^2) \, \|\bar{x}\|^2 \,  \sin^2( \theta )}{n -  \alpha^2 +(\|\bar{x}\| - \alpha  )^2 \sin^2( \theta )} \; . 
    \end{equation} 
\end{lemma}
\noindent {\bf Proof:}  Keep in mind that $ P_{ L^{\perp} }( e ) = P_{L^{\perp}}(\bar{x}) $, because $ \bar{x} \in e + L $.  \vspace{1mm}

 The boundary of $ K_e( \alpha ) $ lies in $ \{ x: (e^T x)^2 - \alpha^2 \, \|x\|^2 = 0 \} \; , $ from which it is seen that the tangent space at $ \bar{x} $ for $ M := (e + L) \cap \partial K_e( \alpha ) $ is 
\begin{align}
   T_{\bar{x}} & = \{ v \in L: (e^T\bar{x}) \, e^T v - \alpha^2 \, \bar{x}^T v = 0 \} \nonumber \\
     &  = \{ v \in L: e^T v = \smfrac{\alpha }{\|\bar{x}\|} \bar{x}^T v \} \quad \textrm{(using $ e^T \bar{x} = \alpha \, \|\bar{x}\| $)} \; .   \label{eqn.ec} \end{align}  
Let $ u $ denote the projection of $ \bar{x} $ onto $ T_{\bar{x}} $.  By (\ref{eqn.ec}), the projection of $ e $ onto $ T_{\bar{x}} $ is $ \smfrac{\alpha }{\|\bar{x}\|} \, u \; , $ from which it follows that the orthogonal complement in $ L $ for $ T_{\bar{x}} $ is the line $ \{ tw: t \in \mathbb{R} \}  $, where
\[      w := (e-\bar{x}) + (1 - \smfrac{ \alpha }{\|\bar{x}\|})u \; , \]
a vector whose length-squared is  $  \|e-\bar{x}\|^2 - (1- \smfrac{ \alpha}{\|\bar{x}\|})^2 \, \|u\|^2  $.  \vspace{1mm}

Again using $ e^T \bar{x} = \alpha \, \|\bar{x}\| $, the length-squared of the projection of $ \bar{x} $ onto $ \{ tw: t \in \mathbb{R} \}  $ is seen to be 
\[  \left(  \frac{w^T \bar{x}}{\|w\|} \right)^2 = \frac{\left(   \alpha \|\bar{x} \| - \|\bar{x}\|^2 + (1 - \smfrac{ \alpha }{\|\bar{x}\|}) \|u\|^2 \right)^2 }{\|e-\bar{x}\|^2 - (1- \smfrac{ \alpha }{\|\bar{x}\|})^2 \|u\|^2} \; . \]
As the length-squared of the projection of $ \bar{x} $ onto $ L $ is this quantity plus $ \|u\|^2 $, the projection of $ \bar{x} $ onto $ L^{\perp} $  thus satisfies 
\begin{align*}
\| P_{ L^{\perp} }(\bar{x}) \|^2 & = \|\bar{x}\|^2 - \|u\|^2 - \frac{ \left(   \alpha \, \|\bar{x} \| - \|\bar{x}\|^2 + (1 - \smfrac{ \alpha }{\|\bar{x}\|}) \|u\|^2 \right)^2 }{\|e-\bar{x}\|^2 - (1- \smfrac{ \alpha }{\|\bar{x}\|})^2 \|u\|^2} \\
& = \|\bar{x}\|^2 \, \sin^2( \theta ) - \frac{ \left(   \alpha \, \|\bar{x} \| - \|\bar{x}\|^2 + (1 - \smfrac{ \alpha }{\|\bar{x}\|}) \|\bar{x}\|^2 \cos^2( \theta ) \right)^2 }{ n - 2 \alpha \|\bar{x}\| + \|\bar{x}\|^2 - (1- \smfrac{ \alpha }{\|\bar{x}\|})^2 \|\bar{x}\|^2 \cos^2( \theta )} \\
& = \|\bar{x}\|^2 \, \sin^2( \theta ) - \frac{ \left(   (\alpha \, \|\bar{x} \| - \|\bar{x}\|^2) (1 - \cos^2( \theta ) \right)^2 }{ n - \alpha^2 + (\|\bar{x}\| - \alpha )^2 (1 - \cos^2( \theta ))} \\
& = \|\bar{x}\|^2 \sin^2(  \theta ) \,  \left(   1 - \frac{( \|\bar{x}\| - \alpha )^2 \sin^2( \theta ) }{n- \alpha^2 + (\|\bar{x}\| - \alpha )^2 \sin^2( \theta ) } \right) \\
& = \frac{(n- \alpha^2) \|\bar{x}\|^2 \sin^2( \theta )}{n- \alpha^2 + (\|\bar{x}\| - \alpha )^2 \sin^2( \theta )}   \; . 
\end{align*} 
Since $ P_{ L^{\perp} }(e) = P_{ L^{\perp} }(\bar{x}) $, the proof is complete.  \hfill $ \Box $
 \vspace{3mm}

  \begin{lemma}  \label{lem.ed}
  Assume $ 0 < \beta  \leq \alpha < \sqrt{n} $, and assume $ e $ satisfies $ \| e \| = \sqrt{n} $.  For every point $ \bar{s}  \in K_e( \beta   )^* $, the cone $ K_e( \alpha   )^* $ contains the ball $ B( \bar{s} ,r) $ where
\[  r = \smfrac{1}{n} \, \| \bar{s} \| \, \left( \alpha  \, \sqrt{  n - \beta^2} - \beta  \, \sqrt{  n - \alpha^2} \,  \right) \; . \]  
\end{lemma}
\noindent {\bf Proof:}  Let $ \delta := \sqrt{ n - \alpha^2} $ and $ \gamma := \sqrt{ n - \beta^2}$.  Thus, $ K_e( \alpha )^* = K_e( \delta ) $ and $ K_e( \beta )^* = K_e( \gamma ) $.  Our goal is to show that if $ \bar{s} \in K_e( \gamma ) $, then $ K_e( \delta ) $ contains the ball $ B( \bar{s}, r) $ where
  \[  r = \smfrac{1}{n} \, \| \bar{s} \| \, \left( \gamma   \, \sqrt{ n - \delta^2} - \delta   \, \sqrt{  n - \gamma^2} \right) \; . \]  

A non-zero vector $ \bar{s}  $  is in $ K_e( \gamma ) $  if and only its angle with $ e $ does not exceed $ \arccos( \gamma / \sqrt{n}) $.  Moreover, for $ \bar{s}  \in K_e( \gamma ) $ ($  \subseteq K_e( \delta ) $), the point in $ \partial K_e( \delta ) $ closest to $ \bar{s}  $ is the projection of $ \bar{s}  $ onto the ray whose angle with $ \bar{s}  $ is smallest among all rays whose angle with $ e $ is $ \arccos( \delta/\sqrt{n} ) $.  Clearly, the angle between $ \bar{s}  $ and this ray is $ \arccos( \delta/ \sqrt{n}) - \theta $, where $ \theta $ is the angle between $ \bar{s}  $ and $ e $.  The distance from $ \bar{s}  $ to the boundary of $ K_e( \delta ) $ thus satisfies
\begin{align*}
\textrm{distance} & = \| \bar{s} \| \,   \sin \Big( \arccos( \delta/ \sqrt{n}) - \theta \Big) \\
& = \| \bar{s} \| \, \Big(  \sin \big( \arccos( \delta/\sqrt{n} \big) \, \cos( \theta ) - \cos \big( \arccos( \delta/ \sqrt{n}) \big) \, \sin( \theta ) \Big) \\
& = \| \bar{s} \| \, \left(  \sqrt{  1 - \delta^2/n } \, \cos ( \theta ) - ( \delta/ \sqrt{n}) \, \sin ( \theta ) \right) \\
& \geq \| \bar{s}  \| \,  \left(  \sqrt{  1 - \delta^2/n } \, ( \gamma/ \sqrt{n} )  - ( \delta/ \sqrt{n}) \, \sqrt{  1 - \gamma^2/n} \right) \\
& = \smfrac{1}{n} \, \| \bar{s} \| \,  \left( \gamma  \sqrt{  n - \delta^2}   - \delta  \, \sqrt{  n - \gamma^2} \right) \; .  
\end{align*}
completing the proof.  \hfill $ \Box $
 \vspace{3mm}

 \noindent {\bf Proof of Proposition~\ref{prop.ea}:} 
    Choose any matrix $ A $ whose nullspace is $ L $.  Define $ b := Ae $, and note $ b \neq 0 $ (because, by assumption, $ e \notin L $). Consider the following primal-dual pair of optimization problems:
 \[  \begin{array}{rl} \min_x & \bar{s}^T x \\
       \textrm{s.t.} & Ax = b \\
                   & x \in K_e( \beta ) \end{array} \qquad  
                   \begin{array}{rl} \max_{y,s} & b^T y \\
                    \textrm{s.t.} & A^T y + s = \bar{s}  \\
                     & s \in K_e( \beta )^* \end{array} \]
By assumption, $ \bar{s} $ lies in the interior of $ K_e( \beta )^* $, and thus the dual problem is strictly feasible.  Since the primal problem is feasible (indeed, $ e $ is a feasible point), the primal problem thus has an optimal solution.  Denote an optimal solution by $ x' $. \vspace{1mm}

Since, $ \bar{s} $ lies in the interior of $ K_e( \beta )^* $ but not in $ L^{\perp} $ (by assumption), $ x' $ is the unique optimal solution and lies in $ \partial K_e( \beta ) $.    Moreover, $ x' \neq 0 $ (because $ b \neq 0 $). Thus it makes sense to speak of the tangent space to $ \partial K_e( \beta ) $ at $ x' $. \vspace{1mm}

Let $ T_{x'} $ denote the tangent space to $ M' := (e + L) \cap \partial K_e( \beta ) $ at $ x' $.  Of course $ T_{x'} \subseteq L $.  Also, $  \dim( T_{x'} ) = \dim(L) - 1 $ (because $ e $ lies in the interior of $ K_e( \beta ) $). \vspace{1mm}

  Optimality of $ x' $ implies that every feasible direction $ v $ that leaves the objective value unchanged must lie in  $ T_{x'} $, that is $ \{ v \in L: \bar{s}^T v = 0 \} \subseteq T_{x'} $.  Thus, because $ \dim(T_{x'}) = \dim(L) - 1 $, we have
  \[     T_{x'} = \{ v \in L: \bar{s}^T v = 0 \} \; , \]
  the important matter being that the set on the right makes no reference to the value $ \beta $. \vspace{1mm}

  In the same manner, then, for $ \bar{x} $ which in the proposition is assumed to be optimal for the primal problem in which $ K_e( \alpha ) $ appears rather than $ K_e( \beta ) $, we have that $ \bar{x} $ is the unique optimal solution, $ 0 \neq \bar{x} \in \partial K_e( \alpha ) $, and  the tangent space $ T_{ \bar{x}} $ at $ \bar{x}$  to  $ M := (e+L) \cap \partial K_e( \alpha) $  satisfies
  \[  T_{ \bar{x}} = \{ v \in L: \bar{s}^T v = 0 \}. \]  
 Thus, $ T_{ \bar{x}} = T_{x'} $.  
\vspace{1mm}

Hence we can infer by Lemma \ref{lem.eb}  that for every $ v \neq 0 $ in the tangent space $ T_{ \bar{x}} $, the angle $ \theta $ between $ v $ and $ \bar{x}  $ satisfies $ | \cos ( \theta ) | \leq \beta/ \alpha $. \vspace{1mm}

Consequently, by Lemma \ref{lem.ec},  
\begin{equation}   \label{eqn.ee}
   \| P_{L^{\perp} }( e ) \| \geq {\mathcal C}_1 := \left(  \frac{(n - \alpha^2) \| \bar{x} \|^2  \Big( 1 - \big(\frac{\beta}{\alpha} \big)^2 \Big)}{n - \alpha^2 + ( \| \bar{x} \| - \alpha )^2 \Big( 1 - \big(\frac{\beta}{\alpha} \big)^2 \Big)} \right)^{1/2} \; .  \end{equation} 
(Here we are using the fact that the expression on the right in (\ref{eqn.ee}) increases as $ r := 1 - \big(\smfrac{\beta}{\alpha} \big)^2 $ increases -- indeed, a simple derivative calculation shows any function of the form $ r \mapsto \frac{ar}{b + cr} $, where $ a, b > 0 $, is strictly increasing everywhere that it is defined.) \vspace{1mm}

On the other hand, by Lemma \ref{lem.ed}, the point 
\[  \bar{s}' := \bar{s} -  {\mathcal C}_2 \, \| \bar{s} \| \, \smfrac{1}{\| P_{ L^{\perp}  }(e) \| } P_{ L^{\perp} }(e) \]
is contained in $ K_e( \alpha )^* $, where
\[   {\mathcal C}_2 := \smfrac{1}{n }  \left( \alpha  \, \sqrt{  n - \beta^2} - \beta  \, \sqrt{  n - \alpha^2} \right)   \; . \]
Clearly, however, 
\begin{align*}
  e^T ( s - \bar{s}') & = {\mathcal C}_2 \, \| \bar{s} \| \, \smfrac{1}{\| P_{ L^{\perp}  }(e) \| } \,  e^T P_{ L^{\perp} }(e)  \\
                & = {\mathcal C}_2 \, \| \bar{s} \| \,  \| P_{ L^{\perp}  }(e) \|  \\
                & \geq  {\mathcal C}_1 \, {\mathcal C}_2 \, \| \bar{s} \| \; ,
                \end{align*}
 thus concluding the proof.  \hfill $ \Box $

\section{{\bf  Proof of the Main Theorem for Semidefinite Programming}} \label{sect.f}

In proving the theorem, we rely on notation that allows the proof to be easily extended to hyperbolic programming in general, the topic of subsequent sections. \vspace{1mm}

For the trace inner product we now write $ \langle \; , \; \rangle  $. Thus, for example, the objective function $ \tr( CX ) $ becomes $ \lin C, X \rin $.  \vspace{1mm}

Recall that for $ E \in \Sympp $, we defined the inner product $ \lin U, V \rin_E := \tr( E^{-1} U E^{-1} V ) \; . $  A different way to express this inner product is 
\[   \lin U, V \rin_E  = \lin U, H(E)[V] \rin \; , \]
where $ H(E) $ is the Hessian at $ E $ for the barrier function 
\[  f(X) = - \ln \det(X) \; , \]
  that is, the linear automorphism of $ \Sym $ given by 
\[  H(E)[V] = E^{-1} V E^{-1}  \]
  (c.f., \cite[Chap.~2.2.1]{renegar2001mathematical}).   \vspace{1mm}

The gradient at $ E $ for the barrier function $ f $ is 
\[  g(E) = -E^{-1} \; . \]
  Thus,
\[ 
   \tr( E^{-1} X) =  - \lin g(E), X \rin =  \lin E, X \rin_E  \; . 
\] 
Moreover,  for $ E \in \swath( \alpha) $, the optimal solutions $ X_E = X_E( \alpha) $ and $ (y_E, S_E) = (y_E( \alpha), S_E( \alpha)) $ are related according to 
\begin{align} 
  S_E  & = \smfrac{\tr(C, E - X_E) }{n - \alpha^2} \left( E^{-1} - \smfrac{\alpha^2}{\tr(E^{-1} X_E) } E^{-1}   X_E  E^{-1} \right) \quad \textrm{(by (\ref{eqn.cc}))} \nonumber \\
&  = \smfrac{ \lin C, E - X_E \rin  }{n - \alpha^2} \left( -g(E) + \smfrac{\alpha^2}{\lin g(E), X_E \rin} H(E)[X_E]   \right) \nonumber \\
& =  \smfrac{ \lin C, E - X_E \rin  }{n - \alpha^2} \, H(E)\left[ E - \smfrac{\alpha^2}{\lin E, X_E \rin_E} X_E \right] \; . \label{eqn.fb} 
\end{align}
\vspace{2mm}

Now we begin the proof of the \hyperlink{targ_main_thm_sdp}{Main Theorem for Semidefinite Programming}.  \vspace{1mm} 

As in the statement of the Main Theorem, assume $ 0 < \alpha < 1 $ and $ E \in \swath( \alpha ) $.  Letting $ X_E = X_{E}( \alpha ) $ and $ (y_E,S_E) = (y_E( \alpha), S_{E}( \alpha )) $, define $ E' $ according to  
\begin{equation}  \label{eqn.fc}
   E' = \smfrac{1}{1 + t_E} \, \left( E + t_E X_E \right) \; , 
   \end{equation} 
where $ t_E = t_{E}( \alpha ) $ minimizes the convex quadratic polynomial 
\begin{equation}  \label{eqn.fca}
   q(t) = \tr \left(   \, \left(  ( E + t X_E) \, S_E \right)^2 \, \right) \; . 
   \end{equation} 
An immediate implication of Corollary \ref{cor.db}  is that $ E' \in \swath( \alpha  ) $, thus establishing the first claim of the Main Theorem.   \vspace{1mm}

For the second claim (i.e., primal objective monotonicity), observe that because $ t_E > 0 $ (Corollary~\ref{cor.db}), $ E' $ has strictly better objective value than $ E $ simply because $ X_E $ does. \vspace{1mm}

For the third claim (i.e., dual objective monotonicity), again we call upon Corollary~\ref{cor.db}, which implies $ S_E $ lies in $ K_{E'}( \alpha )^* $ . Hence, $ (y_E, S_E) $ is feasible for $ \qp_{E'}( \alpha )^* $, implying $ b^T y_E \leq b^T y_{E'} \; , $ where $ y_{E'} = y_{E'}( \alpha) \; . $   Dual objective monotonicity is thus established.  (In fact, the monotonicity is strict, a proof of which begins by noting Corollary~\ref{cor.db}  shows $ S_E $  lies in the {\em interior} of  $ K_{E'}( \alpha)^* $.) \vspace{1mm}

For the final claim of the Main Theorem, as there let $ E_0, E_1, \ldots $, be the sequence obtained by recursively applying (\ref{eqn.fc}), beginning at $ E_0 \in \swath( \alpha) $.  Of course the sequence is well-defined, as we have already established that if $ E_i $ lies in $  \swath( \alpha ) $, then so does $ E_{i+1} $.  Our goal is to show for every for $ i = 0, 1, \ldots  \; , $  
\begin{equation}  \label{eqn.fd} 
      \frac{ \mathrm{gap}_{E_{j+1}}}{ \mathrm{gap}_{E_{j}}} \leq 1 - \frac{\kappa}{\kappa + \sqrt{n}} \;     \quad \textrm{for $ j = i $ or $ j = i+1 $ (possibly both)} \; ,   
      \end{equation} 
      where
      \[ \kappa := \alpha \sqrt{\smfrac{1- \alpha }{8} } \; . \]  
\vspace{1mm}

To simplify notation, let $ X_i = X_{E_i}( \alpha ) $, $ y_i = y_{E_i}(\alpha) $, $ S_i = S_{E_i}( \alpha ) $, and use $ t_i $ to denote the minimizer of the strictly-convex quadratic polynomial 
 \begin{equation}  \label{eqn.fda}
  q_i(t) = \tr( \, (( E_i + t X_i) S_i)^2 \,) \; . 
  \end{equation} 
Thus, 
\[  E_{i+1} = \smfrac{1}{1+t_i} ( E_i + t_i X_i) \; . \]
\vspace{1mm}

We make use of the simple relation
\[  \textrm{gap}_{i+1} = \textrm{gap}_i - \lin C , E_i - E_{i+1} \rin  - b^T (y_{i+1} - y_i) \; . \]
In particular, since the primal and dual objective values are monotonic, 
\begin{equation}  \label{eqn.fe}
            \frac{\textrm{gap}_{i+1}}{\textrm{gap}_i} \leq 1 - \frac{ \lin  C,  E_i - E_{i+1} \rin }{ \textrm{gap}_i} 
    \end{equation}
and 
\begin{equation}  \label{eqn.ff}
            \frac{\textrm{gap}_{i+1}}{\textrm{gap}_i} \leq 1 - \frac{ b^T ( y_{i+1} - y_i) }{ \textrm{gap}_i} \; . \end{equation}
            
To establish (\ref{eqn.fd}), we can fix $ i = 0 $, because $ E_0 $ is an arbitrary point in $ \swath( \alpha) $. We consider two cases: 
\begin{enumerate}

\item Either $ t_0 > \kappa/ \sqrt{n} $ or $ t_1  > \kappa/\sqrt{n} \; . $    
\item Both $ t_0 \leq \kappa/ \sqrt{n} $ and $ t_1 \leq \kappa/\sqrt{n} \; . $
\end{enumerate} 
 \vspace{1mm}

For the first case, proving the duality gap reduction (\ref{eqn.fd}) is easy.  Indeed,
  \begin{align*}
  \frac{\mathrm{gap}_{i+1}}{\mathrm{gap}_i} & \leq    
 1 -  \frac{ \lin C,  E_i - E_{i+1} \rin }{\textrm{gap}_i} \quad \textrm{(by (\ref{eqn.fe}))}  \\
 & = 1 - \frac{\lin C, E_i - E_{i+1} \rin }{\lin C, E_i - X_i \rin  } \quad \textrm{(by (\ref{eqn.cg}))} \\ 
                                              & = 1 - \frac{t_i}{1 + t_i}  \qquad \textrm{(since $ E_{i+1} = \smfrac{1}{1+t_i} \, ( E_i + t_i X_i) $)} 
                                              \\
& < 1 - \frac{\kappa}{\kappa + \sqrt{n}} \quad \textrm{(assuming $ t_i > \kappa/\sqrt{n} $)}         
                                 \end{align*}
 \vspace{1mm}
 
To complete the proof of the Main Theorem, it now only remains to show
\begin{equation} \label{eqn.ffffz}
 t_0, t_1 \leq \kappa/ \sqrt{n} \quad \Rightarrow \quad  \frac{ \mathrm{gap}_{E_1}}{ \mathrm{gap}_{E_{0}}} \leq  1 - \frac{\kappa}{\kappa + \sqrt{n}}  \; . 
 \end{equation}
 Henceforth, assume $ t_0, t_1 \leq \kappa/ \sqrt{n} $. \vspace{1mm}

Let $ {\mathcal L}  $ be the nullspace of $ A $, and let $ {\mathcal L}^{\perp} $ be its orthogonal complement with respect to the inner product $ \langle \; , \; \rangle \; . $  We will also have occasion to refer to the orthogonal complement of $ {\mathcal L}  $ with respect to the inner product $ \langle \; , \; \rangle_{E_1} $, for which we use the notation $ {\mathcal L}^{{\perp}_{E_1}} $. It is easily verified that
\[     {\mathcal L}^{{\perp}_{E_1}} =  \{ H(E_1)^{-1}[V]: V \in {\mathcal L}^{\perp} \} \; .  \]
    
    Using $ E_1 + {\mathcal L} $ as an abbreviation for the set $ \{ E_1 + V: V \in {\mathcal L} \} $, the optimization problem $ \qp_{E_1}( \alpha ) $ can be written
\[ \left. \begin{array}{rl}
        \min & \lin C, X \rin  \\
    \textrm{s.t.} & X \in E_1 + {\mathcal L} \\
          & X \in K_{E_1}( \alpha )  
          \end{array} \right\} \quad \qp_{E_1}( \alpha ) \; . \]
Observe that $ X_1 $, the optimal solution, also is optimal for
\begin{equation}  \label{eqn.fh}
  \begin{array}{rl}
        \min & \lin S_0, X \rin  \\
    \textrm{s.t.} & X \in E_1 + {\mathcal L} \\
          & X \in K_{E_1}( \alpha ) \; . 
          \end{array} \end{equation} 
Indeed, this problem differs from $ \qp_{E_1}( \alpha ) $ only in having $ S_0 $ in the objective rather than $ C $.  However, as 
\begin{equation}  \label{eqn.fi}
  C - S_0 \in {\mathcal L}^{\perp} \quad  \textrm{(because $ A^* y_0 + S_0 = C $)} \; , 
  \end{equation}  
  the two objectives give the same ordering to feasible points (in fact, $ \lin C , X - X' \rin  = \lin S_0, X - X' \rin  $ for all feasible $ X, X' $). Thus, $ X_1 $ in indeed  \vspace{1mm} optimal for (\ref{eqn.fh}). 

Let $ \bar{S}_0 := H(E_1)^{-1}[S_0] \; . $   Trivially, $ \lin S_0, X \rin  = \lin \bar{S}_0, X \rin_{E_1} $ for all $ X \in \Sym $.  Consequently, the optimization problem (\ref{eqn.fh})  can be expressed 
\begin{equation}  \label{eqn.fj}
  \begin{array}{rl}
        \min & \lin \bar{S}_0, X \rin_{E_1}  \\
    \textrm{s.t.} & X \in E_1 + {\mathcal L} \\
          & X \in K_{E_1}( \alpha ) \; . 
          \end{array} \end{equation} 
 \vspace{1mm}

In a different vein, observe Corollary \ref{cor.db}  with $ E = E_0 $ shows    
\[ S_0 \in \int(K_{ E_1}( \beta  )^*)  \quad \textrm{for }   \beta = \alpha \sqrt{  \smfrac{1 + \alpha}{2}} \; . \]
It follows that  $  \bar{S}_0 \in  \int(K_{E_1}( \beta  )^{*_{E_1}}) \; ,  $ 
where $ K_{E_1}( \beta  )^{*_{E_1}} $ is  the cone dual to $ K_{E_1}( \beta  ) $ with respect to the inner product $ \langle \; , \; \rangle_{E_1} $. Indeed, 
\[  K_{E_1}( \beta  )^{*_{E_1}} = \{ H(E_1)^{-1}[S]: S \in K_{E_1}( \beta  )^* \}  \; ,  \]
  a simple consequence of the definition of dual cones and the definition of the inner product $\langle \; , \; \rangle_{E_1} \; .  $  \vspace{1mm}
  
We now make use of Proposition~\ref{prop.ea} which applies to Euclidean spaces generally.  In particular, as $ X_1 $ is optimal for  (\ref{eqn.fj}), and since $ \bar{S}_0 \in  \int(K_{E_1}( \beta  )^{*_{E_1}})   $ (by the immediately preceding paragraph), we can apply Proposition~\ref{prop.ea}  with the inner product $ \langle \; , \; \rangle_{E_1} $, and with  the substitutions $ \bar{s} = \bar{S}_0 $, $ \bar{x} = X_1 $, $ e = E_1 $ and $ L = {\mathcal L} $.  We thereby find:
\begin{quote}
 There exists $ \bar{S}' \in (\bar{S}_0 + {\mathcal L}^{{\perp}_{E_1}}) \cap K_{E_1}( \alpha )^{*_{E_1}} $ satisfying 
\begin{gather}  
 \lin E_1, \bar{S}_0 - \bar{S}' \rin_{E_1} \geq {\mathcal C}_1 \, {\mathcal C}_2   \, \| \bar{S}_0 \|_{E_1} \label{eqn.fk} \\
\textrm{where} \quad   
      {\mathcal C}_1  =  \left(  \frac{(n - \alpha^2) \| X_1 \|_{E_1}^2  \Big( 1 - \big(\frac{\beta}{\alpha} \big)^2 \Big)}{n - \alpha^2 + ( \| X_1 \|_{E_1} - \alpha )^2 \Big( 1 - \big(\frac{\beta}{\alpha} \big)^2 \Big)} \right)^{1/2} \label{eqn.fl}  \\ \textrm{and} \quad  
{\mathcal C}_2  := \smfrac{1}{n }  \left( \alpha  \, \sqrt{  n - \beta^2} - \beta  \, \sqrt{  n - \alpha^2} \right) \; .  \qquad  \label{eqn.fm} 
\end{gather}  
\end{quote}

  Let $ S' := H(E_1)[\bar{S}'] $.  From  $ \bar{S}' \in \bar{S}_0 + {\mathcal L}^{\perp_{E_1}} $ follows $ S' \in S_0 + {\mathcal L}^{\perp } $, and hence $ S' \in C + {\mathcal L}^{ \perp } $ (by (\ref{eqn.fi})) -- thus, there exists $ y' $ for which $ A^* y' + S' = C $. Moreover, from $ \bar{S}' \in K_{E_1}( \alpha )^{*_{E_1}} $ follows $ S' \in K_{E_1}( \alpha )^* $, and thus $ (y',S') $ is feasible for $ \qp_{E_1}( \alpha )^* $. \vspace{1mm}

Consequently, $ (y_1, S_1) $, the optimal solution for $ \qp_{E_1}( \alpha)^* $, satisfies
\begin{align*}
  b^T y_1 - b^T y_0 & \geq b^T y' - b^T y_0 \\
                   & = A( E_1)^T (y' - y_0) \\
                   & = \lin  E_1 ,  A^*( y' -  y_0) \rin  \\
                   & = \lin  E_1 ,    S_0 - S' \rin  \\
                   & = \lin E_1, \bar{S}_0 - \bar{S}' \rin_{E_1} \\
                   & \geq {\mathcal C}_1 \, {\mathcal C}_2   \, \| \bar{S}_0 \|_{E_1} \quad \textrm{(by (\ref{eqn.fk}))} \\
                   & \geq {\mathcal C}_1 \, {\mathcal C}_2   \, \smfrac{1}{\sqrt{n} } \, \lin E_1 , \bar{S}_0 \rin_{E_1} \quad \textrm{(by Cauchy-Schwarz for $ \langle \; , \; \rangle_{E_1} $)} \\
& = {\mathcal C}_1 \, {\mathcal C}_2   \, \smfrac{1}{\sqrt{n} } \, \lin  E_1 , S_0 \rin  \quad \textrm{(using $ \bar{S}_0 = H(E_1)^{-1}[S_0]$ )}   \; .        
 \end{align*}
However, 
 \begin{align*}
\lin  E_1,  S_0 \rin     &  = \smfrac{1}{1 + t_0} \, \lin  E_0 + t_0 X_0, S_0 \rin     \\
                    & =  \smfrac{1}{1 + t_0} \, \lin  E_0,  S_0 \rin     \qquad \textrm{(by complementarity)} \\
    & = \smfrac{1}{1 + t_0} \, \smfrac{\lin C, E_0 - X_0 \rin   }{n - \alpha^2} \, \lin    E_0 \,  , \, E_0 - \smfrac{\alpha^2}{\lin E_0, X_0 \rin_{E_0} }  X_0  \rin_{E_0}    \qquad \textrm{(by (\ref{eqn.fb}))} \\
& = \smfrac{ \lin C, E_0 - X_0 \rin    }{1 + t_0}   \\
& = \smfrac{\textrm{gap}_0 }{1 + t_0} \qquad \textrm{(by (\ref{eqn.cg}))} \\  
& \geq \smfrac{\sqrt{n} }{ \kappa + \sqrt{n} } \, \textrm{gap}_0 \qquad \textrm{(because $ t_0 \leq \kappa/ \sqrt{n} $, by assumption)} \; .         
   \end{align*}
Hence   
 \[   b^T y_1 - b^T y_0 \geq \frac{ {\mathcal C}_1 \, {\mathcal C}_2   \,  \textrm{gap}_0 }{ \kappa + \sqrt{n} } \; .  \]
 Consequently, by (\ref{eqn.ff}),
\[ 
    \frac{\textrm{gap}_1}{ \textrm{gap}_0}  \leq 1 - \frac{ {\mathcal C}_1 \, {\mathcal C}_2  }{ \kappa + \sqrt{n}} \; \; .  
\] 
To complete the proof of (\ref{eqn.ffffz}) (and complete the proof of the Main Theorem), it thus suffices to show
\begin{equation}  \label{eqn.fna}
  {\mathcal C}_1 \, {\mathcal C}_2  \leq \kappa \; . 
 \end{equation}

According to Corollary \ref{cor.db}, $ t_1 \geq \smfrac{1}{2} \, \alpha/ \| X_1 \|_{E_1} \; . $ Hence, since $ t_1 \leq \kappa/\sqrt{n} $ (by assumption), we have 
\begin{equation}  \label{eqn.fr} 
  \|X_1\|_{E_1} \,  \geq \, \frac{\alpha }{ 2 t_1}  \, \geq \, \frac{ \alpha \sqrt{n}}{ 2 \kappa } \, = \, \sqrt{\frac{2n}{1- \alpha}}   \; . 
\end{equation} \vspace{1mm}

A trivial consequence of (\ref{eqn.fr}) and $ 0 < \alpha < 1 $ (as assumed in the theorem), is  
\[  (\|X_1\|_{E_1} - \alpha)^2  \leq \|X_1\|_{E_1}^2 \; .  \]
By this, along with $ n - 1 <  n - \alpha^2 < n $, is immediately seen from (\ref{eqn.fl})  that
\begin{align*}
  {\mathcal C}_1 & >  \left(  \frac{(n - 1) \| X_1 \|_{E_1}^2  \Big( 1 - \big(\frac{\beta}{\alpha} \big)^2 \Big)}{n  + \| X_1 \|_{E_1}^2 \Big( 1 - \big(\frac{\beta}{\alpha} \big)^2 \Big)} \right)^{1/2}  \\
& =  \left(  \frac{(n - 1) \| X_1 \|_{E_1}^2  ( 1 - \alpha) }{2n  + \| X_1 \|_{E_1}^2 ( 1 - \alpha )} \right)^{1/2} \quad \textrm{(substituting $ \beta = \alpha \sqrt{\smfrac{1 + \alpha }{2} } $ )}      \\
& \geq \sqrt{\frac{ n-1}{2}} \qquad  \textrm{(by (\ref{eqn.fr}))} \; .    
  \end{align*}
On the other hand, according to (\ref{eqn.fm}),
\begin{align*}
 {\mathcal C_2} & :=  \smfrac{1}{n }  \left( \alpha  \, \sqrt{  n - \beta^2} - \beta  \, \sqrt{  n - \alpha^2} \right)  \\
& > \smfrac{1}{n} (\alpha - \beta ) \sqrt{n - \beta^2} \qquad  \textrm{(because $ 0 < \beta < \alpha < 1$)} \nonumber \\
& >  \smfrac{\sqrt{n-1} }{n} \, (\alpha - \beta )     \\
& = \smfrac{\sqrt{ n-1}}{n} \, \alpha \, \left(   1 - \sqrt{\smfrac{1 + \alpha }{2}} \right) \; . 
 \end{align*}
Thus,
\begin{align*}
  {\mathcal C}_1 \, {\mathcal C}_2  & > \smfrac{n-1}{n \sqrt{2} } \,  \alpha \,  \left( 1 - \sqrt{\smfrac{1+ \alpha }{2} } \right) \\
   & \geq \smfrac{1}{2 \sqrt{2} } \,  \alpha  \, \left( 1 - \sqrt{\smfrac{1+ \alpha }{2} } \right) \; , \\
   & > \smfrac{1}{2 \sqrt{2} } \,  \alpha  \, \sqrt{1 - \alpha}   \\
   & = \kappa \; ,  
   \end{align*}
  where the last inequality can be verified by substituting $ \gamma = \sqrt{ \smfrac{1 + \alpha }{2}} $ and observing 
  \begin{gather}   \sqrt{1 - \alpha} < 1 - \sqrt{ \smfrac{1+ \alpha }{2}} \quad \textrm{for $ 0 < \alpha < 1 $} \nonumber \\
\Leftrightarrow \nonumber \\
   1 - (2 \gamma^2 - 1)  < ( 1 - \gamma)^2 \quad \textrm{for $ 1/ \sqrt{2} < \gamma < 1 $} \; , 
  \label{eqn.fs}  \end{gather}   
   and then verifying the validity of (\ref{eqn.fs}) via the quadratic formula.  Having established (\ref{eqn.fna}), the proof of the Main Theorem is now complete.   \hfill $ \Box $

\section{{\bf Hyperbolic Programming: Introduction}} \label{sect.g}

We now shift attention to the setting of hyperbolic programming, which includes semidefinite programming as a special case.  As will become clear, the algorithm naturally extends to general hyperbolic programs.  The complexity analysis also immediately extends, with the exception of Corollary~\ref{cor.db}, which played a crucial role.  In order to extend Corollary~\ref{cor.db}  to general hyperbolic programs, we resort to using an especially deep result, the Helton-Vinnikov Theorem (discussed in \S~\ref{sect.j}).  \vspace{1mm}

A {\em hyperbolic program} is a convex optimization problem of the form
\begin{equation}  \label{eqn.ga} 
  \left. \begin{array}{rl}
\min & \lin c, x \rin \\ \textrm{s.t.} & Ax =b \\  & x \in \Lambdap \end{array} \right\}  \, \hp 
\end{equation} 
where $ \Lambdap $ is the closure of a ``hyperbolicity cone'' for a ``hyperbolic polynomial.''  \vspace{1mm}

Letting $ {\mathcal E} $ denote a finite-dimensional Euclidean space, a  {\em hyperbolic polynomial} $ p: {\mathcal E} \rightarrow \mathbb{R} $ is a homogeneous (real) multivariate polynomial for which there exists a vector $ e \in {\mathcal E} $ with $ p(e) \neq 0 $ and for which the restriction of $ p $ to each line in direction $ e $ results in a  univariate polynomial with only real roots -- that is, for all $ x \in {\mathcal E} $, the univariate polynomial $ t \mapsto p(x + te) $ has only real zeros.  The polynomial $ p $ is said to be ``hyperbolic in direction $ e $.'' \vspace{1mm}

For example, $ p: \mathbb{R}^n \rightarrow \mathbb{R} $ defined by $ p(x) = x_1 \cdots x_n $ is easily seen to be hyperbolic in direction $ e$ for all points $ e \in \mathbb{R}^n  $ having no coordinates equal to zero.  (Even if $ e $ has some coordinates equal to zero, still the univariate polynomials $ t \mapsto p(x +te) $ have no real roots for all $ x $; however, the critical property $ p(e) \neq 0$ is not \vspace{1mm} satisfied.) 

Another standard example is $ p(x) := x_n^2 - \sum_{j=1}^{n-1} x_j^2 $, which is hyperbolic in all directions $ e  $ for which $ e_n^2 > \sum_{j=1}^{n-1} e_j^2 $. \vspace{1mm}

A paramount example is $ p(X) := \det(X) $ for $ X \in {\mathcal E} = \Sym $, the vector space of $ n \times n $ symmetric matrices.  Since the eigenvalues of a symmetric matrix are all real, and are the roots of its characteristic polynomial, immediate is the fact that $ p $ is hyperbolic in direction $ E = \pm I $. \vspace{1mm}
  
For a polynomial $ p $ hyperbolic in direction $ e $, the connected component of $ \{ p \neq 0 \} $ containing $ e $ is a special set, known as the {\em hyperbolicity cone}, which we denote $ \Lambdapp \; . $  Clearly, $ \Lambdapp $ is indeed a cone (i.e., closed under multiplication by positive scalars), because $ p $ is homogeneous. The cone $ \Lambdap $ in a hyperbolic program (\ref{eqn.ga}) is the closure of a hyperbolicity cone $ \Lambdapp \; . $  \vspace{1mm}

  In the examples above, $ \Lambdap $  is, respectively, the non-negative orthant (assuming all coordinates of $ e $ are positive), the second-order cone, and the definite cone (either positive definite or negative definite, depending on whether $ E = I $ or $ E = -I $). Each of these sets is not only a cone, it is a {\em convex} cone.  Such is the case for {\em all} hyperbolicity cones.   \vspace{1mm}

More than 50 years ago, G\aa rding \cite{garding1959inequality}  showed that if $ p $ is hyperbolic in direction $ e $, then $ p $ is hyperbolic in direction $ e' $ for every $ e' \in \Lambdapp \; , $ and he derived various consequences, including that $ \Lambdapp $ is convex. Thus, {\em  every} hyperbolic program (\ref{eqn.ga}) is a convex optimization problem.  \vspace{1mm}

It was G\"{u}ler \cite{guler1997hyperbolic}  who brought G\aa rding's results to the attention of the optimization community, and initiated the study of \vspace{1mm} hyperbolic programming.   
 
The above examples of hyperbolicity cones are particularly special in that each one is symmetric, or ``self-scaled'' in the terminology of Nesterov and Todd \cite{nesterov1997self}.  For these examples, primal-dual algorithms apply, including the primal-dual affine-scaling methods discussed in \S\ref{sect.a}.  In general, however, hyperbolicity cones are not even self-dual, let alone self-scaled, and thus hyperbolicity cones are beyond the reach of primal-dual methods, except in some cases where it is known how to ``lift'' the cone into the realm of self-scaled cones \cite{branden2014hyperbolicity, saunderson2012polynomial}, albeit typically at the cost of introducing a not insignificant number of additional variables.  \vspace{1mm}

In any case, even if by increasing the number of variables, some hyperbolic programs can be lifted to the setting of self-scaled cones, doing so is contrary to one of the most promising features of hyperbolic programming, whereby an optimization problem -- say, a semidefinite program -- can be ``relaxed'' to a hyperbolic program whose cone is for a lower degree polynomial (see \cite{renegar2006hyperbolic}  for details).  Solving the relaxed optimization problem potentially provides an approximate solution to the original problem, the semidefinite program.  Moreover, with the relaxed problem having the hyperbolicity cone for a lower degree polynomial, our affine-scaling algorithm is able to solve it faster -- at least in theory -- than primal-dual methods are able to solve the semidefinite program. (An interesting and open  research direction is to explore whether in this way, the limited size of general semidefinite programs solvable in practice can be notably increased in \vspace{1mm} situations where approximate solutions suffice.) 

In the next section, we develop technical preliminaries in a particularly elementary manner that allows us to avoid making reference to the interior-point literature.  In \S\ref{sect.i}, we present the natural generalization of the algorithm to hyperbolic programming, and we state the Main Theorem for Hyperbolic Programming.  Section~\ref{sect.j} discusses the Helton-Vinnikov Theorem, and uses it to generalize Corollary~\ref{cor.db}. In \S\ref{sect.k}, the generalized corollary is shown to be quickly combined with results from sections \ref{sect.d}  and \ref{sect.f} to establish the Main Theorem for Hyperbolic Programming, of which a special case is the Main Theorem for Semidefinite Programming already proven.  Section~\ref{sect.k} closes with a conjecture. 

\section{{\bf Hyperbolic Programming: Technical Preliminaries}}  \label{sect.h} 

We now lay simple groundwork for presenting the algorithm for general hyperbolic programming, and for analyzing the algorithm.  The preliminaries we develop below are special cases of standard results in the interior-point method literature, but the structure provided by hyperbolic polynomials is the basis for proofs that are particularly easy.  As these proofs have not been recorded in the literature, we present them. \vspace{1mm}

Henceforth, fix $ p $ to be a hyperbolic polynomial of (total) degree $ n $, and $ \Lambdapp $ to be its hyperbolicity cone in the Euclidean space $ {\mathcal E} $ having inner product $ \langle \; , \; \rangle \; . $  Let $ \Lambdap $ denote the closure of $ \Lambdapp \; . $  \vspace{1mm}

We rely on the key result of G\aa rding \cite{garding1959inequality}  for which elementary proofs have appeared various places in the literature (e.g., \cite{guler1997hyperbolic, renegar2006hyperbolic}) and hence no proof is given here. \vspace{1mm}

\hypertarget{targ_garding_thm}{}
\begin{garding_thm}  The polynomial $ p $ is hyperbolic in direction $ e $ \\
$ \textrm{~} $ \quad \qquad  \qquad  \qquad  \qquad  \qquad   for {\em every} $ e \in \Lambdapp \; . $  Moreover, $ \Lambdapp $ is convex.
\end{garding_thm}
\addtocounter{prop}{1}

Fix $ f $ to be the function
\[   f(x) := - \ln p(x) \; , \]
with domain $ \Lambdapp \; . $ Let $ g(e) $ denote the gradient of $ f $ at $ e \in \Lambdapp \; , $ and let $ H(e) $ be the Hessian.\footnote{Keep in mind that the gradient depends on the particular inner product, as is clear upon recalling $ g(e) $ is the vector satisfying 
\[ 0 = \lim_{ \| \Delta e \| \downarrow 0} \frac{f( e + \Delta e) - f(e) - \lin g(e), \Delta e \rin }{\| \Delta e \|} \; . \]
Since $ g(e) $ depends on the inner product, so does the Hessian, as it is defined to be the linear operator $ H(e): {\mathcal E} \rightarrow {\mathcal E} $ for which
\[ 0 = \lim_{ \| \Delta e \| \downarrow 0} \frac{\| g(e + \Delta e) - g(e) - H(e)[ \Delta e] \|}{\| \Delta e \|}  \; .  \]}
 \vspace{1mm}

We always assume the closed cone $ \Lambdap $ to be ``pointed,'' i.e., contains no nontrivial subspaces -- equivalently, the dual cone $ \Lambdap^* $ has nonempty interior, where
\[ \Lambdap^* := \{ s \in {\mathcal E} : \lin x, s \rin \textrm{ for all $ x \in \Lambdap $} \} \; . \]
    Thus,  $ \Lambdap $ is ``regular,'' meaning both it and its dual cone have nonempty interior.    The assumption of regularity is standard in the ipm literature (and can be removed at the expense of increasing the number of secondary details). \vspace{1mm}
    
    The following elementary lemma together with G\aa rding's Theorem yields all results in this section. \vspace{1mm}

\begin{lemma}  \label{lem.hb}
Assume $ e \in \Lambdapp $ and assume $ L $ is a two-dimensional subspace of $ {\mathcal E} $ containing $ e $.  There exists a linear transformation $ M: L \rightarrow \mathbb{R}^n $ such that $ Me = {\bf 1} $ (vector of all ones), and for all $ x \in L $, 
\begin{equation}  \label{eqn.ha}
   p(x) = p(e) \, \prod_{i=1}^n y_i \quad  \textrm{where y = Mx} \; . \end{equation} 
\end{lemma}
\noindent {\bf Proof:}  Think of $ p|_L $ ($ p $ restricted to $ L $) as being a polynomial in two variables -- in other words, give $ L $ a coordinate system.  Since $ p|_L $ is homogeneous, the set $ \{x \in \mathbb{R}^{2}:  p|_L(x)  = 0 \} $ is either the origin or a union of real lines through the origin, that is, either the origin or a union of real hyperplanes $ H_i = \{ x \in \mathbb{R}^{2}: \alpha_i^T x = 0 \} $, where $ 0 \neq \alpha_i \in \mathbb{R}^{2} $.  Again due to homogeneity, the value of $ p|_L $ is zero everywhere on the complex hyperplane $ \{ x \in \mathbb{C}^2: \alpha_i^T x = 0 \} $.  Since linear polynomials are (trivially) irreducible, it follows that $ x \mapsto \alpha_i^T x $ is a factor of $ p|_L $.  We claim that all factors of $ p|_L $ arise in this way, i.e., arise from real hyperplanes $ H_i $. \vspace{1mm}

Indeed, assume otherwise, that is, assume $ p|_L(z) = q(z) \prod_{j=1}^m \alpha_i^T z $ where $ m < n $ and where $ q $ is a non-constant polynomial which has no linear factors $ x \mapsto  \beta^T x $ ($ 0 \neq \beta \in \mathbb{R}^{2} $). Since $ p|_L $ is homogeneous, clearly so is $ q $. 
Thus, by the reasoning of the previous paragraph, $ \{ q = 0 \} $ consists only of the origin. But then, for any $ x \in L $ which is not a scalar multiple of $ e $, the univariate polynomial  $ t \mapsto p|_L(x + te) $ has fewer than $ n$  real roots (counting multiplicities), contradicting that $ p $ is hyperbolic in direction $ e $ (G\aa rding's Theorem). \vspace{1mm}

We have shown $ p|_L(x) = \prod_{i=1}^n \alpha_i^T x $ for some vectors $ 0 \neq \alpha_i \in \mathbb{R}^{2}  $. To conclude the proof, simply let $ M := \diag(Ne)^{-1} N \; ,  $ where $ N $ is the $ n \times 2 $ matrix with $ i^{th} $ row $ \alpha_i^T $.    \hfill $ \Box $ \vspace{3mm}

We now derive various consequences of Lemma~\ref{lem.hb}  and \hyperlink{targ_garding_thm}{G\aa rding's Theorem}.  
\vspace{3mm}

\begin{lemma} \label{lem.hc} If $ M $ is the linear transformation in the conclusion to Lemma~\ref{lem.hb}, then
\[  \Lambdapp \cap L \, = \, \{ x \in L: Mx \in \mathbb{R}^n_{\pplus} \} \; . \] 
\end{lemma} 
\noindent {\bf Proof:} Recall that $ \Lambdapp $ is the connected component of $ \{ x \in {\mathcal E}: p(x) \neq 0 \} $ containing $ e $.  Since $ \Lambdapp $ is convex, it follows that $ \Lambdapp \cap L $ is the connected component of $ \{ x \in L: p(x) \neq 0 \} $ containing $ e $.  With this characterization of $ \Lambdapp \cap L $, the identity is immediate from (\ref{eqn.ha}). \hfill $ \Box $
\vspace{3mm}

\begin{lemma} \label{lem.hd}  If $ M $ is as in the conclusion to Lemma~\ref{lem.hb}, then
$ M $ is one-to-one.
\end{lemma} 
\noindent {\bf Proof:}  If $ M $ was not injective, then $ p $ would be constant on a line through $ e $, due to (\ref{eqn.ha}). The line would be contained in $ \Lambdapp $ (because $ \Lambdapp $ is the connected component of $ \{ p \neq 0 \} $  containing $ e $). This would contradict our standing assumption that $ \Lambdap $ is pointed.  \hfill $ \Box $
\vspace{3mm}

 \begin{lemma} \label{lem.he} 
 If $ M $ is as in the conclusion to Lemma~\ref{lem.hb}, then for all $ u, v \in L \; , $
 \[  \lin g(e), v \rin = - \mathbf{1}^T M v \quad \textrm{and} \quad \lin u, H(e) \,  v \rin = (M u)^T M v \; , \]
 where $ {\bf 1} $ is the vector of all ones. 
\end{lemma}
\noindent {\bf Proof:}      
Let $ F: \mathbb{R}^n_{\pplus} \rightarrow \mathbb{R} $ denote the function $ F(y) = - \sum_j \ln y_j \; . $ Thus,  $ f(x) = F(Mx) - \ln p(e) $. \vspace{1mm}

 Let $ \nabla F $ and $ \nabla^2 F $ be the gradient and Hessian of $ F $ with respect to the standard inner product on $ \mathbb{R}^n $.  Trivially, $ \nabla F( {\bf 1}) = - {\bf 1} $ and $ H( {\bf 1}) = I $, the identity.  
Recalling that $ g $ and $ H $  denote the gradient and Hessian of $ f $ with respect to the inner product $ \langle \; , \; \rangle \; , $ and  letting $ M^* $ denote the adjoint\footnote{The adjoint of $ M $ is the linear operator $ M^*: \mathbb{R}^n  \rightarrow L $ satisfying $ \lin M^* u, v \rin = u^T Mv  $ for all $ u \in \mathbb{R}^n  $, $ v \in L $.} of $ M $,  it holds for all $ v \in L $ that 
\begin{align*}
     \lin g(e), v \rin  & = \lin M^* \, \nabla F(Me)  \, , v \rin  \\
             & =  - \lin M^*{\bf 1}  \, , v \rin   \\
             & = - {\bf 1}^T  (Mv)  \; . 
             \end{align*}
Likewise, for all $ u, v \in L $,             
 \begin{align*}
          \lin u, H(e) v \rin  & = \lin u, M^*  \, \nabla^2 F(Me)[Mv] \, \rin \\
                          & = \lin u, M^*  \,  Mv \,  \rin  \\    
                          & = (Mu)^T (Mv) \; ,
                     \end{align*}
completing the proof.  \hfill $ \Box $
\vspace{3mm}

\begin{prop}  \label{prop.hf} 
The Hessian $ H(e) $ is positive definite for all $ e \in \Lambdapp \; . $ 
\end{prop} 
\noindent {\bf Proof:}  Fix $ e \in \Lambdapp $.  The goal is to show for an arbitrary vector $ 0 \neq  v \in {\mathcal E} $ that $ \lin v, H(e) \,  v \rin > 0 $. However, letting $ L $ be the subspace spanned by $ e $ and $ v $, and letting $ M $ be the resulting linear transformation in Lemma~\ref{lem.hb}, we have from Lemma~\ref{lem.he}   that 
\[    \lin v, H(e) \, v \rin \, = \, \| Mv \|^2 \; . \]
Since $ M $ is one-to-one (Lemma~\ref{lem.hd}), the proof is complete.  \hfill $ \Box $
\vspace{3mm}

Proposition~\ref{prop.hf}  implies, of course, that $ f $ is strictly convex. \vspace{1mm}
  
 Define the ``local inner product'' at $ e \in \Lambdapp $ by 
 \[  \lin u, v \rin_e \, := \, \lin u, H(e) \, v \rin \quad \textrm{for $ u, v \in {\mathcal E} $ } \; . \]
(By Proposition~\ref{prop.hf}, this is indeed an inner product on $ {\mathcal E} $.)  The local inner products $ \langle \; , \; \rangle_e $  are natural for performing complexity analysis of interior-point methods.  These inner products are easily shown to be independent of the initial (arbitrary) inner product $ \langle \; , \; \rangle $ -- that is, the local inner products are dependent only on $ f  $ -- and thus are sometimes referred to as the ``intrinsic inner products'' (c.f., \cite[Chap.~2.1]{renegar2001mathematical}). \vspace{1mm}

\begin{lemma} \label{lem.hg} If $ e \in \Lambdapp \; , $ then $ H(e) \, e = -g(e) \; . $  
\end{lemma}
\noindent {\bf Proof:} Fix $ e \in \Lambdapp \; . $ To prove the lemma, it suffices to show for all $ v \in {\mathcal E} $ that $ \lin v, H(e) \,  e \rin = - \lin v, g(e) \rin \; . $ \vspace{1mm}

Fix $ v \in {\mathcal E} \; . $  Letting $ L $ be the subspace spanned by $ e $ and $ v $, and letting $ M $ be as in the conclusion to Lemma~\ref{lem.hb},
\begin{align*} 
     \lin v, H(e) \, e \rin &  = (Mv)^T M e   \quad \textrm{(by Lemma~\ref{lem.he})} \\
                      & = (Mv)^T \mathbf{1} \\
                      & = - \lin g(e), v \rin \quad \textrm{(by Lemma~\ref{lem.he})} \; , 
\end{align*} 
thus concluding the proof.  \hfill $ \Box $
 \vspace{3mm}

  Let $ \| \; \; \|_e $ be the norm associated with $ \langle \; , \; \rangle_e \; , $ i.e., $ \|v\|_e = \lin v, v \rin_e^{1/2} $.  For $ x \in {\mathcal E} $ and $ r > 0 $, let  $ B_e(x,r) := \{ x' \in {\mathcal E}: \| x' - x \|_e < r \} \; . $  \vspace{1mm}

 The gist of the next lemma is that no point $ e $ in the cone $  \Lambdapp $ considers itself as being ``close'' to the boundary of the cone. \vspace{1mm}

\begin{prop}  \label{prop.hh} 
If $ e \in \Lambdapp \; , $ then $ \|e\|_e = \sqrt{n} $ and  $ B_e(e,1) \subset \Lambdapp \; . $ 
\end{prop} 
\noindent {\bf Proof:}  Assume $ x \in B_e(e,1) $, let $ L $ be the subspace spanned by $ x $ and $ e $, and let $ M $ be as in Lemma~\ref{lem.hb}.  Since $ Me = {\bf 1} $, Lemma~\ref{lem.he}  implies $ \|e\|_e^2 = {\bf 1}^T {\bf 1} = n \; . $ \vspace{1mm}

To prove the containment, it suffices by Lemma~\ref{lem.hc}  to show $ Mx \in \mathbb{R}^n_{\plus} $. Thus, as  Lemma~\ref{lem.he}  implies $ \|x - e\|_e = \| Mx - {\bf 1} \|  \; , $ it suffices to observe that the Euclidean unit ball centered at $ {\bf 1} $ is contained in the strictly-positive orthant. \hfill $ \Box $  \vspace{3mm}

The objects of primary interest for us are not the $ \| \;\; \|_e $-balls, but their close relatives, the closed quadratic cones defined for positive values $ \alpha $ by
\begin{align}
    K_e( \alpha ) & := \{ x \in {\mathcal E}  : \lin e, x \rin_e \geq \alpha \, \|x\|_e \} \label{eqn.hb} \\
                  & \, = \{ x \in {\mathcal E}: \lin e, x \rin_e \geq 0 \, \textrm{ and } \, \lin e, x \rin_e^2 - \alpha^2 \| x \|_e^2 \geq 0 \,  \}  \label{eqn.hc}  \\
& = \{ x \in {\mathcal E} : \angle_e (e,x) \leq \arccos\left(\alpha/ \sqrt{n} \, \right) \} \; ,   \label{eqn.hd} 
                  \end{align}     
where $ \angle_e( e,x) $ denotes the angle between $ e $ and $ x $ measured in terms of the local inner product $ \langle \; , \; \rangle_e \; , $ and where in the final equality we have used $ \| e \|_e = \sqrt{n} $ (by Proposition~\ref{prop.hh}).   If $ \alpha > \sqrt{n} $, then $ K_e( \alpha) = \emptyset \; . $     \vspace{1mm}

The only essential difference between our algorithm and the original affine scaling algorithm due to Dikin is that we replace his use of balls $ B_e(e,r) $ -- known as ``Dikin ellipsoids -- by a use of cones $ K_e( \alpha) $. \vspace{1mm}

\begin{prop}  \label{prop.hi} If $ e \in \Lambdapp \; , $ then \, $ K_e( \sqrt{n-1} \,  ) \subseteq  \Lambdap \subseteq K_e(1 ) \; . $ 
\end{prop} 
\noindent {\bf Proof:} Assume $ 0 \neq x \in K_e( \sqrt{n-1} \, ) $, that is, assume 
\[    \lin e, x \rin_e \geq \sqrt{n-1} \, \| x \|_e \; . \]
The $ \langle \; , \; \rangle_e $-projection of $ e $ onto $ x $ is 
\[  \bar{x} = \smfrac{\lin e, x \rin_e}{\|x\|_e^2} \,  x \; , \]
a positive-scalar multiple of $ x $.  To show $ x \in \Lambdap \; , $  it suffices to show $ \bar{x} \in \Lambdap \; , $ for which in turn it suffices, by Proposition~\ref{prop.hh}, to show $ \| \bar{x} - e \|_e \leq 1 \; . $ However,   
\begin{align*}
\| \bar{x} - e \|_e^2 & =  \| \bar{x} \|_e^2 - 2 \lin e, \bar{x} \rin_e + n \qquad  \textrm{(using $ \|e\|_e = \sqrt{n} $)} \\
& = - \left(  \frac{\lin e, x \rin_e}{\|x\|_e} \right)^2 + n \\
& \leq 1 \; , 
\end{align*} 
completing the proof that $ K_e( \sqrt{n-1} \, ) \subseteq \Lambdap \; . $     \vspace{1mm}

To establish $ \Lambdap \subseteq K_e(1) \; , $  fix $ x \in \Lambdap \; . $   Let $ L $ be the subspace spanned by $ e $ and $ x $, and let $ M $ be as in Lemma~\ref{lem.hb}.  From Lemma~\ref{lem.he}, 
\[   \lin e, x \rin_e = {\bf 1}^T Mx \quad \textrm{and} \quad \| x \|_e = \| M x \| \; . \]  Thus, to show $ x \in K_e(1) $, it suffices to show  $ y = Mx $ satisfies
\begin{equation}  \label{eqn.he}
     {\bf 1}^T y  \geq \|y\|  \; .  
\end{equation} 
However, this inequality  is satisfied by all vectors $ y \in \mathbb{R}^n_{\plus} $ -- indeed, the set consisting of all vectors $ y $ satisfying (\ref{eqn.he})  is convex, and vectors $ y $ on the positive coordinate axes satisfy (\ref{eqn.he}) with equality.  Since $ Mx \in \mathbb{R}^n_{\plus} $ (from Lemma~\ref{lem.hc}), we thus have $ x \in K_e(1) \; ,  $ completing the proof that $ \Lambdap \subseteq K_e(1) \; . $   \hfill $ \Box $  \vspace{3mm} 

The first half of the preceding proof is easily adapted to make precise the relation between Dikin ellipsoids $ B_e(e,r) $ and the quadratic cones -- specifically, for $ r \leq \sqrt{n} \; , $  the smallest cone containing $ B_e(e,r) $ is $ K_e(\sqrt{n - r^2} \, ) \; . $ (When $ r > \sqrt{n} \; , $ the smallest cone is all of $ {\mathcal E} $.) \vspace{1mm}

Clearly, if $ \alpha_1 \leq \alpha_2 $, then $ K_e( \alpha_1 ) \supseteq  K_e( \alpha_2 ) \; . $ Thus, Proposition~\ref{prop.hi} implies for all $ 0 \leq \alpha \leq 1 $ that $ K_e( \alpha ) $ is a ``relaxation'' of $ \Lambdap $ (i.e., contains $ \Lambdap $).  Moreover, $ \alpha = 1 $ is the largest value for which, in general, $ K_e( \alpha ) $ contains $ \Lambdap $ (indeed, consider the hyperbolic polynomial $ p(x) = x_1 \cdots x_n $ whose cone is the positive orthant). Likewise, $ \alpha = \sqrt{n-1} $ is the smallest value which for all hyperbolic polynomials of degree $ n $,  $ K_e( \alpha ) \subseteq \Lambdap \; . $  \vspace{1mm}

In the following sections, we always assume $ 0 < \alpha \leq 1 $, in which case, by the proposition, $ \Lambdap \subseteq K_e( \alpha) $.  An important consequence is that the dual cones satisfy the reverse inclusion,  $ \Lambdap^* \supseteq K_e( \alpha)^* $.  \vspace{1mm}

 We make use of the following characterization of the dual cone $ K_e( \alpha )^* $. \vspace{1mm}

\begin{prop}  \label{prop.hj} 
If $ e \in \Lambdapp $ and $ 0 <  \alpha <  \sqrt{n} \; , $ then 
\[  
    K_e( \alpha)^*  = \{ H(e) \, s: s \in K_e \big( \sqrt{n - \alpha^2} \, \big) \} \; . \]  
\end{prop}
\noindent {\bf Proof:}  
According to (\ref{eqn.hd}),  
\[   K_e( \alpha) = \{ x: \angle_e(e,x) \leq \arccos( \alpha / \sqrt{n}) \} \; . \]
Thus, the dual cone with respect to the inner product $ \langle \; , \; \rangle_e $ is
\begin{align*}
   K_e( \alpha)^{*_e} & :=  \{ s: \lin x, s \rin_e \geq 0 \textrm{ for all $ x \in K_e( \alpha) $} \} \\
& = \{ s: \angle_e(e,s) \leq \arcsin( \alpha / \sqrt{n}) \}   \\
& = \{ s: \angle_e(e,s) \leq \arccos( \sqrt{n - \alpha^2}/\sqrt{n} \, ) \\
& = K_e( \sqrt{n - \alpha^2} \, ) \; , 
\end{align*}
where the last equality is due to  (\ref{eqn.hd}) applied with $ \sqrt{n - \alpha^2} $ substituted for $ \alpha $.    
However, as the local inner product is defined by $ \lin u, v \rin_e = \lin u, H(e) \, v \rin $, it is immediate that the dual cone with respect to the original inner product satisfies
 $ K_e( \alpha)^* = \{ H(e)s: s \in K_e( \alpha)^{*_e} \} \; . $   \hfill $ \Box $
 \vspace{3mm}
 
 We close the section with a brief discussion of the ``eigenvalues of $ x $ in direction $ e $,''  these being the roots of the univariate polynomial 
 \begin{equation}  \label{eqn.hf}
  \lambda \mapsto p(\lambda e - x ) \; . 
  \end{equation}
  Clearly, this generalizes the usual notion of eigenvalues, for which the relevant hyperbolic polynomial is $ X \mapsto \det(X) $ and the direction is $ I $, the identity matrix. \vspace{1mm}

  Let $ \lambda_e(x) \in \mathbb{R}^n $ be the vector of eigenvalues, where the number of times an eigenvalue appears is equal to its multiplicity as a root of the polynomial (\ref{eqn.hf}), and where, for definiteness, the eigenvalues are ordered, say, from smallest to largest.    \vspace{1mm}
 
For reasons that will become clear, determining the step length in our algorithm amounts to minimizing a convex quadratic polynomial $ t \mapsto a t^2 + b t + c $ where $ a $, $ b $ and $ c $ depend on current points $ e $ and $ x $ through the evaluation of the following four symmetric functions at the coordinates of  $ \lambda_e(x) = ( \lambda_1, \ldots, \lambda_n) $:
\[  \sum_j \lambda_j \; , \quad \sum_j \lambda_j^2 \; , \quad \sum_j \lambda_j^3 \quad \textrm{and} \quad \sum_j \lambda_j^4 \; . \]
We observe that computing these values can be done efficiently.  Indeed, first compute the five leading coefficients of the univariate polynomial
\[  t \mapsto  p(x + te) \, = \, \sum_{i=0}^n a_i t^i \, = \, p(e) \, \prod_{j=1}^n (t + \lambda_j) \]
(equivalently, by homogeneity, the five trailing coefficients of the reverse polynomial $ s \mapsto p(sx + e) = \sum_{i=0}^n a_i s^{n-i} $ \, ).
Then
\[       \frac{a_{n-1}}{a_n} = \sum_i  \lambda_i \; , \quad \frac{a_{n-2}}{2a_n} =  \sum_{i < j} \lambda_i \lambda_j \; , \quad \frac{a_{n-3}}{6a_n}  =   \sum_{i < j < k} \lambda_i \lambda_j \lambda_k \; , \quad  \frac{a_{n-4}}{24 \, a_n} =  \sum_{i < j < k < \ell } \lambda_i \lambda_j \lambda_k \lambda_{\ell} \; , \]
from which follow
\[ \sum_j \lambda_j = \frac{a_{n-1}}{a_n} \; , \quad \sum_j \lambda_j^2 = \left( \frac{a_{n-1}}{a_n} \right)^2 - \frac{a_{n-2}}{a_n} \; , \quad \sum_j \lambda_j^3 = \left( \frac{a_{n-1}}{a_n} \right)^3 - \frac{3}{2} \, \frac{a_{n-1}}{a_n} \, \frac{a_{n-2}}{a_n} + \frac{1}{2} \, \frac{a_{n-3}}{a_n} \]
and
\[ \sum_j \lambda_j^4 = \left( \frac{a_{n-1}}{a_n} \right)^4 - 2 \, \left( \frac{a_{n-1}}{a_n} \right)^2 \frac{a_{n-2}}{a_n} + \frac{1 }{2} \, \left( \frac{a_{n-2}}{a_n} \right)^2 + \frac{2}{3} \, \frac{a_{n-1}}{a_n} \, \frac{a_{n-3}}{a_n} - \frac{1}{6} \frac{a_{n-4}}{a_n} \; ,  \]
the first four \hypertarget{targ_newton_girard_identities}{}``Newton-Girard identities.''

\section{{\bf Duality, the Algorithm,  \\
 and the Main Theorem for Hyperbolic Programming}} \label{sect.i} 
 
Consider a hyperbolic program
\[ \left. \begin{array}{rl}
 \min & \lin c, x \rin \\
\textrm{s.t.} & Ax = b \\
& x \in \Lambdap  \end{array}  \right\} \, \hp   \]
where $ \Lambdap $ is a regular cone whose interior is the hyperbolicity cone for a hyperbolic polynomial $ p: {\mathcal E} \rightarrow \mathbb{R} $ of degree $ n $.   We assume $ A $ is of full rank and $ b \neq 0 $ (in particular, the origin is infeasible), and we assume $ c $ does not lie in the image of the adjoint of $ A $ (otherwise all feasible points for $ \hp $ are optimal).   \vspace{1mm}

For $ e \in \Lambdapp $ and $ 0 < \alpha \leq 1 $, the following quadratic-cone optimization problem is a relaxation of $ \hp $ (by Proposition~\ref{prop.hi}): 
\[ \left. \begin{array}{rl}
 \min & \lin c, x \rin \\
\textrm{s.t.} & Ax = b \\
& x \in K_e( \alpha)  \end{array}  \right\} \, \qp_e(\alpha)   \]
\vspace{1mm}

We give the following generalization of Definition~\ref{def.bab}:
\begin{definition} 
\[  \swath(\alpha ) := \{ e \in \Lambdapp: Ae = b \textrm{ and } \qp_e(\alpha) \textrm{ has an optimal solution} \} \] 
\end{definition}  
\vspace{2mm}

For $ \alpha > 0 $ and $ e \in \swath( \alpha) $, let $ x_e = x_e(\alpha) $ denote the unique optimal solution of $ \qp_e( \alpha ) $ (unique because $ K_e( \alpha ) $ is a regular quadratic cone, and because the origin is infeasible). Analogous to the setting of semidefinite programming, our algorithm takes a step from $ e $ towards $ x_e $. \vspace{1mm}

When $ e \in \swath(\alpha) $, computation of $ x_e $ is readily accomplished is a way essentially identical to the one for computing $ X_E $ in the context of semidefinite programming, which relied on the first-order optimality conditions (\ref{eqn.bb}). Now, since
\[   \lin e, x \rin_e^2 - \alpha^2 \|x\|_e^2 = \, \lin g(e), x \rin^2 - \alpha^2  \lin x, H(e) x \rin \quad \textrm{(using Lemma~\ref{lem.hg})} \; , \]
necessary conditions for a point $ x $ to be optimal for $ \qp_e( \alpha) $ are as follows:
\begin{align}
  & 0 = \lin g(e), x \rin^2  - \alpha^2 \, \lin x, H(e) x \rin  \nonumber \\   & Ax = b \nonumber  \\
   & 0 = \lambda c + A^* y +  \lin g(e), x \rin \, g(e)  - \alpha^2 H(e) x \label{eqn.ia}  \\ &
\qquad  \qquad    \textrm{for some $ y \in \mathbb{R}^{m} $ and $ \lambda \in \mathbb{R} $} \; .  \nonumber \end{align}
 Excluding the first equation gives a system of linear equations in the variables $ x $, $ y $ and $ \lambda $ whose solutions form a one-dimensional set.  Parameterizing this one-dimensional set linearly, substituting for $ x $ in the first equation and then using the quadratic formula to determine the roots, results in two candidates $ x' $, $ x'' $, one of which is optimal for $ \qp_e( \alpha ) $ (assuming, as we are, that $ e \in \swath(\alpha) $).  Finally, checking feasibility and objective values of $ x' $ and $ x''$ reveals which is the optimal solution. \vspace{1mm}

The optimization problem dual to $ \qp_e( \alpha) $ is
\[ \left. \begin{array}{rl} \max_{y,s} & b^T y \\
   \textrm{s.t.} & A^* y + s = c \\
           & s \in K_e( \alpha)^* \end{array} \right\}  \qp_e( \alpha)^*  \]
where $ A^* $ is the adjoint of $ A $, and $ K_e( \alpha)^* $ is the dual cone for $ K_e( \alpha) $.  Just as was described for semidefinite programming immediately following (\ref{eqn.ca}), the first-order optimality conditions (\ref{eqn.ia}) provide us with the optimal solution $ (y_e, s_e) = (y_e( \alpha), s_e( \alpha) ) $ of $ \qp_e( \alpha)^* $, namely, $   y_e = - \frac{1}{\lambda } y $ and 
\[ 
  s_e = \smfrac{\lin c, e - x_e \rin }{n - \alpha^2} \, \left( - g(e) + \smfrac{\alpha^2}{\lin g(e),  x_e \rin} H(e) \, x_e \right)  \; , \] 
where $ y $ and $ \lambda $ are the values which together with $ x_e $ satisfy the first-order conditions.  (The generalization of the expression for $ \lambda $ given in (\ref{eqn.caa})  is easily verified to be $ \lambda =  \frac{(n - \alpha^2) \, \lin g(e),  x_e \rin}{\lin c, e - x_e \rin} \; . $) \vspace{1mm}

Lemma~\ref{lem.hg} shows $ H(e)\, e  = - g(e) \; .  $ Consequently, 
\begin{equation}  \label{eqn.ib}
 s_e(\alpha) = \smfrac{\lin c, e - x_e \rin }{n - \alpha^2} \, H(e) \, \left(   e - \smfrac{\alpha^2}{\lin e, x_e \rin_e} \,  x_e \right)   \; , 
  \end{equation} 
  which is the expression for $ s_e = s_e( \alpha) $ on which we will rely.  \vspace{1mm}

Assuming $ 0 < \alpha \leq 1 $, the cone $ K_e( \alpha) $ is a relaxation of $ \Lambdap $, and hence $ K_e( \alpha )^* \subseteq \Lambdap^* $.  Consequently, $ (y_e, s_e) = (y_e( \alpha), s_e( \alpha)) $ is feasible for the optimization problem dual to the hyperbolic program $ \hp $:
\[   \left. \begin{array}{rl} \max & b^T y \\
    \textrm{s.t.} & A^* y + s = c \\
    & s \in \Lambdap^* \end{array} \right\} \hp^* \]
Thus, the duality gap between $ (y_e, s_e) $ and $ e $ (which is feasible for $ \hp $) is of interest.  Exactly as we obtained the identity (\ref{eqn.cg})  in the context of semidefinite programming, now we have
\begin{align*}   
\textrm{gap}_e & := b^T y_e - \lin c, e \rin \nonumber \\
& =  \lin c, e - x_e \rin  
\end{align*} 
(thus,  the difference in the primal objective values for the points $ e $ and $ x_e $ is a duality gap in disguise). \vspace{1mm}

\hypertarget{targ_main_thm_hp}{}
\begin{main_thm_hyperbolic}  
Assume $ 0 < \alpha < 1 $ and $ e \in \swath( \alpha) $.
 Let $ \lambda = \lambda_e(x_e(\alpha)) $ be the vector of eigenvalues of $ x_e( \alpha) $ in direction $ e $, and let $ \tilde{q}  $ denote the strictly-convex quadratic polynomial 
 \begin{equation}  \label{eqn.ic} 
 \tilde{q}(t) = a t^2 + b t + c 
 \end{equation}
  where
\begin{gather}  
a = \Big( \sum_j \lambda_j \Big)^2 \sum_j \lambda_j^2 - 2 \alpha^2 \Big( \sum \lambda_j \Big) \sum_j \lambda_j^3 + \alpha^4 \sum_j \lambda_j^4 \; , \label{eqn.id} \\
 b = 2 \alpha^4 \sum_j \lambda_j^3 - 2 \, \Big( \sum_j \lambda_j \big)^3  \quad \textrm{and} \quad c = (n - \alpha^2) \, \Big( \sum_j \lambda_j \Big)^2 \; . \label{eqn.ie} 
\end{gather} 
Define 
\begin{equation}  \label{eqn.if} 
      e'  = \smfrac{1}{1 + t_e( \alpha)} \big( e + t_e(\alpha ) \,  x_e(\alpha ) \big) \; , 
      \end{equation} 
where $ t_e( \alpha) $ is the minimizer of $ \tilde{q}  $. \vspace{1mm}

Then
\begin{itemize}

\item $ e' \in \swath( \alpha ) $ 
\item $ \lin c, e \rin > \lin c, e' \rin $ \quad (primal objective monotonicity) 
\item $ b^T y_e( \alpha ) \leq b^T y_{e'}(\alpha) $ \quad (dual objective monotonicity)
\end{itemize}
Moreover, if beginning with $ e_0 \in \swath( \alpha) $, the identity (\ref{eqn.if}) is recursively applied to create a sequence $ e_0, e_1, \ldots \; , $ then for every $ i = 0, 1, \ldots \; , $  
\[ 
     \frac{ \mathrm{gap}_{e_{j+1}}}{ \mathrm{gap}_{e_{j}}} \leq  1 - \frac{\kappa}{\kappa + \sqrt{n}}     \quad \textrm{for $ j = i $ or $ j = i+1 $ (possibly both)} \; ,  \]  
 where
 \[ \kappa := \alpha \sqrt{\smfrac{1- \alpha }{8} } \; . \]           
\end{main_thm_hyperbolic}
\vspace{2mm}

The quadratic polynomial $ \tilde{q}  $ is efficiently computed by making use of the \hyperlink{targ_newton_girard_identities}{Newton-Girard identities} presented at the end of \S\ref{sect.h}.  \vspace{1mm}

The proof of the Main Theorem for Hyperbolic Programming is virtually identical to the proof of the theorem in the special case of semidefinite programming, as we explain in \S\ref{sect.k}.  The one major difference is that we have been unable to directly prove the generalization of Corollary~\ref{cor.db}, and thus we resort to making use of a deep result -- the Helton-Vinnikov Theorem --  to accomplish the generalization.

\section{{\bf  An Application of the Helton-Vinnikov Theorem}}  \label{sect.j} 

The Helton-Vinnikov Theorem extends Lemma~\ref{lem.hb}  by allowing the subspace $ L $ to be 3-dimensional.  Whereas Lemma~\ref{lem.hb}  was almost trivial to prove, the Helton-Vinnikov Theorem is extremely deep. \vspace{1mm}

\begin{helton_vinnikov_thm}  
\hypertarget{targ_helton_vinnikov_thm}{} Assume $ e \in \Lambdapp \; , $ assume $ L $ is a 3-dimensional subspace of $ {\mathcal E} $ containing $ e $, and recall $ n $ denotes the degree of $ p $.   There exists a linear transformation $ T: L \rightarrow \Sym $ such that 
 \[     T(e) = I \quad \textrm{and} \quad p(x) = p(e) \, \det(T(x)) \, \textrm{ for all $ x \in L $} \; . \]
\end{helton_vinnikov_thm}

\addtocounter{prop}{1}
 
This is a ``homogeneous version'' of the theorem  proved by Helton and Vinnikov \cite{helton2007linear} (also see \cite{vinnikov2012lmi}).  The homogenous version was first recorded by Lewis, Parrilo and Ramana in \cite{lewis2005lax}, and observed there to settle affirmatively the Lax conjecture (i.e., Lax conjectured that the statement of the above theorem is true). \vspace{1mm}

Br\"{a}nd\'{e}n \cite{branden2011obstructions}  showed that the theorem cannot, in general, be extended to subspaces $ L $ of dimension four.  \vspace{1mm}

Frankly, we would prefer to avoid relying on a result as deep as the Helton-Vinnikov Theorem, but to date, we have been unable to generalize Corollary~\ref{cor.db}  without making use of the theorem.  Following is our generalization of the corollary, to which this section is devoted to proving. \vspace{1mm} 

\begin{cor}  \label{cor.jb} 
Assume $ 0 < \alpha < 1 $ and $ e \in \swath( \alpha) $.  Let $ \tilde{q}  $ denote the quadratic polynomial $ \tilde{q}(t) = a t^2 + b t + c $ whose coefficients are given in (\ref{eqn.id})  and (\ref{eqn.ie}).  Then $ \tilde{q}  $ is strictly convex, and its minimizer $ t_e( \alpha) $ satisfies $ \delta := t_e( \alpha) - \smfrac{1}{2} \alpha / \|x_e( \alpha) \|_e > 0 \; $; moreover,
\[ t_e( \alpha) - \delta \leq t \leq t_e( \alpha) + \delta \quad \Rightarrow \quad e(t) \in \swath( \beta ) \,  \textrm{ and }\,  s_e( \alpha) \in \int( K_{e(t)}( \beta )^* ) \; ,  \]
where $ e(t) = \smfrac{1}{1+t} (e + t x_e( \alpha)) $   and  $ \beta = \alpha \sqrt{ \frac{1 + \alpha}{2}} \; . $ \end{cor}
\vspace{2mm} 

In addition to using the Helton-Vinnikov Theorem in proving the corollary, we rely on two results from \S\ref{sect.d}.  \vspace{1mm}

\begin{prop}  \label{prop.jc} \textnormal{(Proposition \ref{prop.da}  for special case $ E = I $, the identity matrix.)} Assume  $ 0 < \alpha < 1 $ and $ 0 \neq X  \in \partial K_I( \alpha ) $ ($ \subset \Sym $).  Let  \[   S  := \smfrac{1}{n - \alpha^2}  \big( I - \smfrac{\alpha^2}{\tr(X)  }  X   \big)   \; . \]
The quadratic polynomial
\begin{equation} \label{eqn.ja} 
   q(t) =  \tr \left(   \, \left(  ( I + t X ) \, S  \right)^2 \, \right)  
   \end{equation} 
is strictly convex and its  minimizer $ \bar{t}  $ satisfies $ \delta :=   \bar{t} -\smfrac{1}{2} \alpha / \| X  \| > 0 $ (Frobenius norm).  Moreover,
\[  \bar{t} - \delta \leq t \leq \bar{t} + \delta \quad \Rightarrow \quad E(t) \in \Sympp \, \textrm{ and } \, S \in \int\big( K_{E(t)}( \beta )^* \big) \]
where $ E(t) = \smfrac{1}{1+t} (I + tX) $ and $  \beta = \alpha \sqrt{\smfrac{1 + \alpha }{2} } \; . $ 
\end{prop} 
\vspace{3mm}

\begin{prop} \label{prop.jd}  \textnormal{(Corollary \ref{cor.dga}  for the special case $ E = I $.)}  Under the assumptions of Proposition~\ref{prop.jc}, the polynomial $ q $ in (\ref{eqn.ja}) is a positive multiple of 
\begin{equation}  \label{eqn.jb}
     t \mapsto a t^2 + b t + c 
 \end{equation}     
where $ a $, $ b $ and $ c $ are given by the formulas (\ref{eqn.id}) and (\ref{eqn.ie}), but here with $ \lambda = \lambda(X) $ being the vector of eigenvalues for the symmetric matrix $ X $. 
\end{prop}
\vspace{3mm}

\noindent {\bf Proof of Corollary~\ref{cor.jb}:} Assume $ e \in \Lambdapp  $ and assume $ L $ is a 3-dimensional subspace of $ {\mathcal E} $ containing $ e $.  Let $ T: L \rightarrow \Sym $ be a linear transformation as in  the \hyperlink{targ_helton_vinnikov_thm}{Helton-Vinnikov Theorem}; thus, $ T(e) = I $, the identity matrix. We begin by recording a few observations. \vspace{1mm}

Observe first that for each $ x \in L $, the roots of the univariate polynomial $ \lambda \mapsto p( \lambda e - x) $ are the same as the roots of 
\[ 
  \lambda \mapsto \det( T( \lambda e - x)) = \det( \lambda I - X) \quad \textrm{where $ X = T(x) $} \; \; ,  \]
  that is, the eigenvalue vectors are identical, $ \lambda_e(x) = \lambda (X) \; . $   \vspace{1mm}

Observe next that for all $ \bar{e}  \in \Lambdapp \cap L $ and $ u, v \in L $, 
\begin{equation}  \label{eqn.jc} 
   \lin u, v \rin_{ \bar{e}} = \tr( \bar{E}^{-1} U \bar{E}^{-1} V ) \quad \textrm{where $ \bar{E} = T( \bar{e}) $, $ U = T(u) $ and $ V = T(v) $} \; . 
   \end{equation}
Indeed, since by definition, $ \lin u, v \rin_{ \bar{e}} = \lin u, H( \bar{e}) \, v \rin $ where $ H( \bar{e}) $ is the Hessian at $ \bar{e} $ for the function $ f(x) = - \ln p(x) $, the identity (\ref{eqn.jc})  is an immediate consequence of the fact that with respect to the trace inner product on $ \Sym $, the Hessian at $ X $ for the function $ X \mapsto  - \ln \det(X) $ is the linear automorphism $  W \mapsto  X^{-1} W X^{-1} $. \vspace{1mm}

A consequence of (\ref{eqn.jc})  is that for $ \bar{e} \in \Lambdapp \cap L $ and $ x \in L $, 
\begin{equation}  \label{eqn.jd} 
  x \in K_{ \bar{e}}( \gamma ) \quad \Leftrightarrow \quad X \in K_{ \bar{E}}( \gamma ) \; ,  
  \end{equation}
  regardless of the value $ \gamma $. \vspace{1mm}

Assume now that $ 0 \neq x \in \partial K_e( \alpha) $, and hence $ 0 \neq X \in \partial K_I( \alpha) \,  $. Let $ \hat{q} $ be the quadratic polynomial $ \hat{q}(t) = a t^2 + bt + c $ where $ a $, $ b $ and $ c $ are given by the formulas (\ref{eqn.ic}) and (\ref{eqn.id}) using $ \lambda = \lambda_e(x) \; . $      Since $ \lambda_e(x) = \lambda(X) $, $ \hat{q} $ is identical to the polynomial  (\ref{eqn.jb}) in Proposition~\ref{prop.jd}, and thus by that proposition, $ \hat{q} $   is a positive multiple of the polynomial $ q $ specified in Proposition~\ref{prop.jc}.   \vspace{1mm} 

Consequently, Proposition~\ref{prop.jc}  implies $ \hat{q} $ is strictly convex, and that its minimizer $ \hat{t}  $ satisfies $ \hat{\delta} := \hat{t}  - \smfrac{1}{2} \alpha / \|x\|_e > 0 $ (using $ \|x\|_e = \|X\| $, by (\ref{eqn.jc})). Moreover, since $ \Lambdapp \cap L = \{ \bar{e} \in L: T( \bar{e}) \in \Sympp \} $,  from Proposition~\ref{prop.jc}  we infer
\[ 
  \hat{t} - \hat{\delta} \leq t \leq \hat{t} + \hat{\delta} \quad \Rightarrow \quad e(t) \in \Lambdapp \; , 
\] 
where $ e(t) = \smfrac{1}{1+t} ( e + tx) \; . $ \vspace{1mm}

Assume further that $ 0 < \alpha < 1$ and $ e \in \swath( \alpha) $.  We specialize $ x $ to be $ x = x_e = x_e( \alpha) $, the optimal solution of $ \qp_e( \alpha) $, in which case the polynomial $ \hat{q} $ above is precisely the polynomial $ \tilde{q}  $ in the statement of Corollary~\ref{cor.jb}.     In light of the immediately preceding paragraph, to establish Corollary~\ref{cor.jb},  it only remains to show for the dual optimal solution $ (y_e, s_e) = (y_e ( \alpha), s_e( \alpha)) $ that
\begin{equation} \label{eqn.je} 
   t_e  - \delta_e \leq t \leq t_e + \delta_e \quad \Rightarrow \quad  s_e \in \int( K_{e(t)}( \beta )^* ) \; , \end{equation}
where $ e(t) = \smfrac{1}{1+t} (e + t x_e) $ and $ \beta = \alpha \sqrt{\frac{1 + \alpha}{2}} \; . $ (The corollary's claim that $ e(t) \in \swath( \beta) $  -- i.e., the claim that $ \qp_{e(t)}(\beta ) $ has an optimal solution -- follows from strict feasibility of $ (y_e, s_e) $ for the dual problem $ \qp_{e(t)}( \beta)^* $, as implied by (\ref{eqn.je}), and the strict feasibility of $ e(t) $ for the primal problem $ \qp_{e(t)}( \beta ) $.) 
Deserving of mention is that the results of the immediately preceding paragraph can be established with Lemma~\ref{lem.hb}  in place of the Helton-Vinnikov Theorem.  Only in the step below does the power of the theorem become critical.  \vspace{1mm}

For each value $ t $ in the range specified in (\ref{eqn.je}), our goal is to show 
\begin{gather*}  s_e = \smfrac{\lin c, e - x_e \rin }{n - \alpha^2} \, H(e) \,   \left( e - \smfrac{\alpha^2}{\lin e, x_e \rin_e}  \, x_e \right) \\ \textrm{ satisfies } \\\lin z, s_e \rin > 0 \textrm{ for all }  0 \neq z \in K_{e(t)}( \beta ) \; ,  
\end{gather*} 
that is, our goal is to show
\begin{gather}  s := \smfrac{1}{\lin c, e - x_e \rin } H(e)^{-1} s_e = \smfrac{1}{n - \alpha^2} \left( e - \smfrac{\alpha^2}{\lin e, x_e \rin_e} x_e \right)  \nonumber \\
  \textrm{ satisfies } \label{eqn.jf}  \\
  \lin z, s \rin_e > 0 \textrm{ for all }  0 \neq z \in K_{e(t)}( \beta ) \; . \nonumber
 \end{gather}

Fix $ t $ within the specified range, and fix $ 0 \neq z \in K_{e(t)}( \beta ) $.  Let $ L $ be the (\underline{three-dimensional}) subspace spanned by $ e $, $ x_e $ and $ z $.  For $ T $ as in the \hyperlink{targ_helton_vinnikov_thm}{Helton-Vinnikov Theorem}, let $ X_e = T(x_e) $, $ Z = T(z) $ and  \vspace{1mm} $ E(t) = T( e(t)) = \smfrac{1}{1+t} (I + t X_e) \; . $

Observe that $ 0 \neq Z \in K_{E(t)}( \beta  ) $, by (\ref{eqn.jd}). Also observe 
 that the vector $ s $ in (\ref{eqn.jf})  lies in $ L $, and  
\[   S := T(s) = \smfrac{1}{n - \alpha^2} ( I - \smfrac{\alpha^2}{\tr(X_e)} X_e ) \]
(using $ \lin e, x_e \rin_e = \tr( X_e ) $, by (\ref{eqn.jc})  when $ \bar{e} = e $). Thus, according to Proposition~\ref{prop.jc}  (applied with $ X = X_e $), $ S \in \int(K_{E(t)}( \beta )^*) $. It follows, since $ 0 \neq Z \in K_{E(t)}( \beta  ) $, that
\[  \lin z, s \rin_e = \tr( ZS) > 0 \; , \]
thereby establishing (\ref{eqn.jf})  and hence completing the proof.  \hfill $ \Box $
 \vspace{3mm}

\section{{\bf Proof of the Main Theorem for Hyperbolic Programming, \\ and  a Conjecture}}  \label{sect.k} 

The proof of the \hyperlink{targ_main_thm_hp}{Main Theorem for Hyperbolic Programming} is identical to the proof given in \S\ref{sect.f} for the special case of semidefinite programming, except  
\begin{itemize}
\item the barrier function $ f(X) = - \ln \det(X) $ is replaced by the function $ f(x) = - \ln p(x) \; , $   

\item the polynomial $ q(t) = \tr \left(     \, \big((E + t X_E) S_E \big)^2 \,  \right)   $ in (\ref{eqn.fca}) is replaced by the polynomial $ \tilde{q} $ specified in Corollary~\ref{cor.jb} (and similarly for the polynomials $ q_i $ in (\ref{eqn.fda})) 
\end{itemize}
and
\begin{itemize} 
\item every reference to Corollary~\ref{cor.db} is replaced by reference to Corollary~\ref{cor.jb}.
\end{itemize}
Also, of course, the upper case letters used for matrices are replaced by lower case letters; for example, the sequence $ E_0, E_1, E_2, \ldots \,  $ becomes $ e_0, e_1, e_2 , \ldots \; , $ and  
\[    S_E = \smfrac{\lin C, E - X_E \rin }{n - \alpha^2} \, H(E) \left[ E - \smfrac{\alpha^2}{\lin E, X_E \rin_E} X_E \right] \]
becomes
\[ s_e   = \smfrac{\lin c, e - x_e \rin }{n - \alpha^2} \, H(e) \, \left( e - \smfrac{\alpha^2}{\lin e, x_e \rin_e} x_e \right)  \]
(in the context of semidefinite programming, we wrote $ H(E)[V] $ rather than $ H(E) \, V $ to emphasize $ H(E) $ is a linear operator acting on $ V \in \Sym $, not an $ n \times n $ matrix multiplying $ V $). \vspace{1mm}

With these slight changes to the exposition in \S\ref{sect.f}, the proof is complete.   \hfill $ \Box $
 \vspace{5mm}
 
 We close with observations that lead to a conjecture. \vspace{1mm}

   The polynomial $ \tilde{q} $ in the \hyperlink{targ_main_thm_hp}{Main Theorem for Hyperbolic Programming} specializes in the case of semidefinite programming to be a positive multiple of the polynomial $ q $ in the \hyperlink{targ_main_thm_sdp}{Main Theorem for Semidefinite Programming}; specifically,
 \begin{equation}  \label{eqn.ka}
   q(t) = \smfrac{1}{\left(  (n - \alpha^2) \sum_j \lambda_j \right)^2} \, \tilde{q}(t) \; , 
   \end{equation} 
where $  \lambda_1, \ldots, \lambda_n $ are the eigenvalues of the matrix $ E^{-1/2} X_E E^{-1/2} $.  Proposition~\ref{prop.de} shows that for every $ 0 < \beta < 1 $ and $ t > -1 $, 
\begin{equation}  \label{eqn.kb}
   ( E(t) \in \Sympp) \, \wedge \, ( S_E( \alpha) \in \int( K_{E(t)}( \beta )^* ) \quad \Leftrightarrow \quad q(t) < \frac{1}{n - \beta^2} \; . 
   \end{equation} 
Our proofs used only the implication ``$ \Leftarrow $'', but the equivalence ``$ \Leftrightarrow $'' is nice in that it precisely characterizes the largest step $ t > 0 $  allowable if the goal is to ensure that $ (y_E( \alpha), S_E( \alpha)) $ is feasible for the dual quadratic-cone problem $ \qp_{E(t)}( \beta )^* \; . $   \vspace{1mm}

Now, it is tempting to hope that for hyperbolic programming in general, the precise characterization analogous to (\ref{eqn.kb}) is satisfied by the polynomial $ \tilde{q} $, or rather, by the multiple of $ \tilde{q} $ on the right of (\ref{eqn.ka}), where here $ \lambda_1, \ldots, \lambda_n $ are to be interpreted as the eigenvalues of $ x $ in direction $ e $.  However, for this quadratic polynomial, our use of the Helton-Vinnikov Theorem in proving Corollary~\ref{cor.jb} only establishes the implication ``$\Leftarrow$'', not the equivalence ``$ \Leftrightarrow $''.  \vspace{1mm}

The proof of Proposition~\ref{prop.de} readily extends to show that for hyperbolic programming in general, the equivalence ``$ \Leftrightarrow $'' belongs to the function
\begin{equation}  \label{eqn.kc}
   t \mapsto \frac{ \lin s_e( \alpha), H \left(  e + t \, x_e( \alpha) \right)^{-1} s_e( \alpha) \rin}{ \lin e, s_e( \alpha ) \rin^2 }  \; . 
   \end{equation} 
 This function easily is seen to be rational in $ t $, but in general certainly is not quadratic, nor in general does it even extend to an entire function on the reals.  Letting 
 \[   t_{ \max} := \max \{ t: e + t x_e( \alpha) \in \Lambdap \} \; , \]
  we conjecture that the function (\ref{eqn.kc}) is convex on the open interval $ (0, t_{\max} )$ and has a minimizer therein.  However, even if the conjecture is true, there remains the issue of evaluating the function quickly if the function is to play a role in an algorithm which is efficient in practice.  We see no way to accomplish this in general.  We are thus content that the efficiently-computable quadratic polynomial on the right of (\ref{eqn.ka})  at least provides the implication ``$\Leftarrow$'' critical for our proof of the \hyperlink{targ_main_thm_hp}{Main Theorem for Hyperbolic Programming}.

\bibliographystyle{plain}
\bibliography{affine_scaling}

\end{document}